# A hybrid large neighborhood search algorithm for the integrated dial-a-ride problem using electric vehicles

**Yumeng Fang**
Department of Engineering, University of Luxembourg, Esch-sur-Alzette, Luxembourg
Email: yumeng.fang@uni.lu; ORCID: 0000-0002-4137-5920

**Tai-Yu Ma**
Luxembourg Institute of Socio-Economic Research (LISER), 11 Porte des Sciences, 4366 Esch-sur-Alzette, Luxembourg
Email: tai-yu.ma@liser.lu; ORCID: 0000-0001-6900-098X

**Abstract**

Integrating demand-responsive mobility services with transit systems is recognized as a practical and effective strategy to mitigate their impact on traffic congestion and the environment. This study develops an efficient hybrid metaheuristic to solve the integrated dial-a-ride problem by utilizing electric vehicles to minimize operational costs and customer travel time. Customer transfer inconvenience is restricted by a maximum intermodal transfer time to synchronize demand-responsive buses' arrival and transit departures. The proposed metaheuristic addresses the challenges of integrating demand-responsive vehicle routing and charging operations with fixed-route transit systems with capacitated charging stations and partial recharge. We benchmarked our algorithm against a state-of-the-art mixed-integer programming solver on instances with 10-50 customers and two transit lines. Our approach achieves solutions that are, on average, 23.8% better in solution quality within around 2 minutes, outperforming those obtained by the solver using an 8-hour computational time limit. We evaluate the impact of various system parameters to bridge the gap between theory and practice. The results suggest that, from the operator's perspective, while the integrated dial-a-ride service reduces vehicle kilometers traveled, the used fleet size may not necessarily be reduced when ensuring high-quality service for passengers. Moreover, operating the integrated systems is more beneficial in areas with dense transit networks, compared with increases in transit frequency. The findings provide valuable insights for developing integrated dial-a-ride services in practice.

# 1 Introduction

Transport network companies such as Uber, Lyft, and DiDi have significantly contributed to demand-responsive transport (DRT) services, transforming the urban mobility landscape and altering the public's travel habits. While these services complement public transportation and offer increased convenience to passengers, their rapid expansion has also introduced new challenges and problems. Henao and Marshall (2019) conducted a survey in Denver, USA, analyzing passengers who used on-demand services such as Lyft, UberX, Lyft Line, and UberPool. The survey revealed that DRT services increased the vehicle kilometers traveled (VKT) by approximately 83.5% compared with their absence, with most induced VKT stemming from public transport demand. Similarly, Tirachini and Gomez-Lobo (2020) found that DRT trips primarily replace demand from public transport and create induced demand, with each new trip increasing VKT by an average of 1.79 km. The increase in VKT is mainly due to deadhead driving to pick up passengers and a shift in demand from public transportation. Schaller (2021) draws similar conclusions using data from Uber and Lyft in both dense and suburban areas of the US. The increased VKT leads to greenhouse gas emissions and potentially causes traffic congestion. Expanding infrastructure, implementing traffic zone restrictions, or introducing congestion pricing are common strategies for managing the impact of DRT services (Behroozi, 2023). Alternatively, integrating DRT service with existing mass transit offers a sustainable solution, as VKT increases are mitigated only when ride-hailing is used for first-mile/last-mile trips connecting public transportation (Ma et al., 2019; Schaller, 2021).

An Integrated Demand-Responsive Transport (IDRT) service is a flexible transport service that combines/integrates demand-responsive vehicles with fixed-route transit to better accommodate customer travel needs. The IDRT service has been conceptualized since the 1980s. The associated optimization problem, known as the Integrated Dial-a-Ride Problem (IDARP), involves designing vehicle routes and schedules that minimize system operating costs while ensuring a desired level of service (Wilson et al., 1976). Häll et al. (2009) proposed a mixed-integer linear programming (MILP) formulation for the IDARP with high-frequency transit service. This model was later extended to include timetabled transit lines, considering transfer times at transit stations (Posada et al., 2017). Subsequently, Posada and Häll (2020) and Molenbruch et al. (2021) developed efficient metaheuristics to solve this timetabled IDARP. However, the majority of the existing studies focus on minimizing vehicle travel cost from the operators' perspective. Considering the composite objective with both vehicle routing time and passenger travel time brings additional challenges as it needs to satisfy passengers' maximum riding time constraints involving different multimodal service options (e.g., bus + train) (Aldaihani and Dessouky, 2003). Moreover, minimizing bus routing costs might increase passengers' riding time, requiring the development of a tailored algorithm to address this issue (Su et al., 2024).

With the accelerating climate crisis, public transport operators have started electrifying their fleets. This trend brings additional challenges for jointly optimizing charging operations and vehicle routing in the context of IDARP using electric vehicles (named EIDARP). Fang et al. (2024) is the first study to incorporate electric vehicles (EVs) in the IDARP. The authors developed a MILP model for a more generalized EIDARP where the trade-off between operator cost and customer inconvenience is handled in the objective function. Furthermore, customer requests can be rejected with high penalty costs. For charging, they consider a partial recharge policy with charging station capacity constraints. The developed model was tested over a set of small instances with up to 20 customers and compared with the benchmark (Posada et al., 2017). However, developing efficient solution algorithms remains a research challenge for the practical applications of electric integrated demand-responsive transport (EIDRT) systems.

The contributions of this study are summarized as follows. First, we develop a hybrid large neighborhood search (LNS) metaheuristic to efficiently solve the EIDARP. The hybrid LNS introduces several key features tailored to the considered problem: 1) proposing a departure-expanded graph to model the synchronization between demand-responsive vehicles and transit service; 2) hybridizing a set of LNS operators with tailored local search operators and using the deterministic annealing (DA) framework for the worse solution acceptance parameters update to address the complex EIDARP; 3) extending the widely-used eight-step scheme of classical DARP (Parragh, 2011) for evaluating solution feasibility in the context of integrated multimodal service options using EVs; 4) developing efficient charging schedule algorithms using partial recharge policy under capacitated charging station constraints. The hybrid LNS is tested on a series of test instances with up to 100 requests and compared with the MILP solutions obtained from the commercial MILP solver (Gurobi). The results demonstrate that the proposed hybrid LNS outperforms the MILP solver and efficiently obtains high-quality solutions. A series of numerical experiments is conducted to assess the impacts of key EIDARP model parameters, including fleet size, transit service characteristics, charging operations, and system configurations, from both customers' and operators' perspectives. This provides valuable insights for the practical applications of the model.

The remainder of this paper is organized as follows. Section 2 presents the literature review. It follows the problem description and modeling of EIDARP in Section 3. Section 4 introduces the hybrid LNS algorithm. Section 5 presents the test instance generation, algorithm parameter tuning, and computational results. Section 6 analyzes the impacts and performance of EIDARP with respect to a set of key model parameters. Finally, conclusions are drawn, and future extensions are discussed.

## 2      Literature review

We focus our literature review on the following three aspects: i) integrated DRT system modeling and solution approach; ii) solution feasibility evaluation procedures for DARP and IDARP; iii) metaheuristics for solving IDARP and related hybridization techniques. A summarized table comparing related studies is presented at the end of this section.

### 2.1      Integrated demand-responsive transport services

Early studies for IDRT system design integrated demand-responsive transport services into fixed-route transit service as a first-mile/last-mile service. For example, Liaw et al. (1996) formulated the IDRT problem as a bi-modal DARP, focusing on operator cost minimization in which an online and offline decision support system was developed to optimize vehicle routing. Hickman and Blume (2001) considered the trade-off of operator costs and passenger inconvenience in IDRT vehicle route optimization (i.e., maximum travel time and number of transfers). The authors proposed a two-stage heuristic to solve the problem efficiently: (i) determining which customers can use integrated options (demand-responsive vehicles plus fixed-route transit); (ii) scheduling the first- and last-mile legs of trips to minimize vehicle routing costs. A realistic case study was conducted with thousands of passenger requests, showing that adopting IDRT would improve system efficiency and significantly reduce operator costs. Subsequently, Aldaihani and Dessouky (2003) developed a Tabu algorithm to optimize vehicle routing while minimizing the increase in customer travel times. Their algorithm was tested on a paratransit service with up to 155 requests. The results demonstrated a 16.6% reduction in VKT, albeit with an 8.7% increase in customer travel time.

The MILP formulation of IDARP was first brought up by Häll et al. (2009). The authors extended classical DARP with fixed-route transit lines. The service frequency is assumed to be unlimited. Posada et al. (2017) extended the IDARP of Häll et al. (2009) by incorporating a timetabled transit service where passengers can get on subsequent transit services (defined by the timetables) given the demand-responsive vehicles' arrival time at transit stations. As the problem is NP-hard, it was only solved exactly for small instances with up to 6 requests. Later, the authors developed an Adaptive Large Neighborhood Search (ALNS) algorithm to solve the problem with up to 145 requests and one transit line (Posada and Häll, 2020). Molenbruch et al. (2021)

developed an LNS for a similar IDARP where walking to/from nearby transit stations is excluded. The authors introduced an efficient LNS algorithm to optimize vehicle arrivals at transit stations, minimizing customer waiting time while adhering to maximum ride time constraints. The algorithm was tested on a set of IDRP test instances adapted from the DARP instances of Braekers et al. (2014). The results show that the proposed method can obtain good solutions efficiently. The authors comprehensively analyzed system performance under various demand scenarios, service frequencies, transit line speeds, and detour factors. Melis et al. (2024) investigated a demand-responsive bus routing problem to minimize total customer travel times (no maximum detour constraints for each passenger). Their model assumes high-frequency transit services where waiting time at transit stations is neglected. An insertion-based solution approach is proposed, along with a speed-up procedure for determining whether parts of the customer journeys can be replaced by fixed-route transit services. In logistics, the intermodal transport problem is a similar routing problem. It involves two types of transport modes, and vehicle routing and scheduling need to be synchronized to minimize transport costs (Ghilas et al., 2016, 2018).

## 2.2 Feasibility evaluation for DARP and IDARP

One of the computation bottlenecks for solving DARP or IDARP is the solution feasibility evaluation for verifying a series of constraints concerning vehicles (capacity, time windows at pickup and/or drop-off locations, energy (if EVs are involved), etc.) and passengers (maximum ride time, pairing constraints that a passenger can only be picked up /dropped off by the same vehicle). For the DARP, the widely used eight-step scheme (Cordeau and Laporte, 2003; Parragh et al., 2010) can efficiently verify the feasibility of solutions. This algorithm evaluates the feasibility of routes in terms of vehicle capacity, time windows, and customer ride time. If a time-related constraint is violated, a vehicle departure delay process is initiated using forward slack times, which calculates the maximum allowed vehicle departure time delay without violating time window constraints for the remaining part of the route. Then repair its infeasibility by delaying vehicle departure based on this information.

Braekers et al. (2014) pointed out that the eight-step scheme increases computational complexity and suggested a preliminary check on these constraints to reduce unnecessary computations. Cordeau and Laporte (2007) formulate the scheduling problem of a route as a MILP to determine optimal vehicle departure time so time windows are satisfied and excess passenger ride times are minimized. To efficiently solve this subproblem, Gschwind and Drexl (2019) proposed a constant-time feasibility check that evaluates only the pickup and delivery nodes when a request is inserted. The approach achieves an average speedup factor of 3.8 compared with the classical eight-step scheme. Bongiovanni et al. (2024) developed a novel algorithm for the excess ride time minimization of DARP. The proposed approach corrects the incorrect infeasibility declaration of the eight-step scheme and can find an optimal solution for most cases. The authors further extend their method for the electric DARP (EDARP) by incorporating a battery management heuristic for vehicle charging scheduling. This heuristic is based on the idea of inserting charging operations as early as possible and as much as possible. If no feasible charging schedule can be found, the route is energy infeasible. Su et al. (2023) developed an exact solution algorithm for EDARP based on the forward labeling algorithm, with a more generalized objective function that considers the weighted sum of vehicle travel time and passengers' excess ride times. The authors developed a DA algorithm by integrating the exact feasibility evaluation scheme to solve the EDARP efficiently.

For the IDARP, the solution feasibility evaluation presents unique challenges due to the need to minimize transfer times at transfer stations and coordinate vehicle routing when customer journeys involve two different demand-responsive vehicles (first-mile and last-mile connections to transit stations). This will complicate the feasibility evaluation for customers' maximum ride time constraints. Hickman and Blume (2001) addressed this issue by identifying potential transfer points with respect to customers' pickup and drop-off locations, formulating customer journey planning as a shortest path problem with time window constraints. Aldaihani and Dessouky (2003) addressed this transfer synchronization problem by identifying each customer's feasible paths, which involve a transit station pair connecting the customer's origin and destination. The objective is to minimize the total ride time of customers. Posada and Häll (2020) applied the constant-time feasibility check of Gschwind and Drexl (2019), but did not provide the details for the feasibility check of IDARP solutions. Molenbruch et al. (2021) formulated the IDARP feasibility check problem as a simple temporal problem, modelled on a directed graph where the nodes represent service times at different locations and the arcs represent time-related constraints such as time windows, maximum ride time, and transit synchronization. The authors extend the Bellman-Ford algorithm to efficiently verify the feasibility of solutions of IDARP.

As we consider an IDRT service using EVs, charging operations can be scheduled at positions of the vehicle route when there are no onboard passengers. Our idea is similar to that of Bongiovanni et al. (2024), which involves inserting charging operations as early as possible while utilizing forward slack times to meet time window constraints when inserting charging operations (Ma et al., 2024).

### 2.3 Large neighborhood search and its hybridization

LNS, initially introduced by Shaw (1998), has been successfully applied for solving a wide range of vehicle routing problems, where the neighborhoods of a promising solution are defined by a set of "destroy" (remove parts of a solution) operators and of "repair" (reinsert the removed part of a solution) operators. Adaptive LNS (ALNS) further extends this approach by dynamically selecting destroy and repair operators based on a weighting system (Ropke and Pisinger, 2006). ALNS has been widely applied for various routing problems, such as vehicle routing, location routing, PDP, and DARP (Mara et al., 2022). Turkeš et al. (2021) investigated the benefit of the adaptive layer in ALNS (using meta-analysis) on improving the solution quality and found an average improvement of 0.14% compared with LNS.

Different from Variable Neighborhood Search (VNS), which systematically alternates between diversification and intensification during its solution search, (A)LNS primarily leverages diversification. Hybrid approaches that integrate intensification strategies with LNS have been explored in the literature. For instance, Alinaghian and Shokouhi (2018) incorporated improvement phases using VNS within the framework of ALNS to solve the multi-depot multi-compartment vehicle routing problem. Akpunar and Akpinar (2021) hybridized ALNS with VNS to solve the capacitated location routing problem, achieving enhancements in computational time and solution quality compared with the standalone ALNS. Another approach to balance diversification and intensification of LNS is through a dynamic acceptance criterion (e.g., simulated annealing or DA (also known as threshold accepting)). This method allows temporarily accepting worse solutions with a small probability, which is controlled by a temperature parameter. By updating this parameter in the solution search process, diversification and intensification can be adapted to find better solutions. The main advantage of the DA is that the number of algorithmic parameters to tune is relatively small compared with the VNS. Moreover, it performs equally well, but with only a fraction of the computational time required by the simulated annealing approach (Dueck and Scheuer, 1990). As a result, the DA has been widely used to solve the DARP and its variants (Braekers et al., 2014; Ma et al., 2024; Su et al., 2023). Santini et al. (2018) conducted a systematic evaluation of nine acceptance criteria within an adaptive LNS framework and found that DA is one of the best acceptance criteria to use.

(A)LNS has been applied to solving IDRT problems in both freight transportation (Ghilas et al., 2016) and passenger transportation (Molenbruch et al., 2021; Posada and Häll, 2020). These studies introduced specific destroy and repair operators tailored to address specific constraints related to integrating transit routes and demand-responsive vehicles. For instance, Ghilas et al. (2016) and Molenbruch et al. (2021) utilized a similarity indicator between pairs of customers, considering whether the same transit stations are used for the first or last mile by the repair operator. Posada and Häll (2020) proposed two IDARP-specific operators: i) repair operator examining the possibility of introducing a fixed-route transit leg into a customer trip; ii) destroy operator removing integrated customer routes in a descending order of route costs. These studies tested their algorithms on real-world test instances and benchmarked them with instances of DARP or PDP. However, none of these studies incorporates tailored local search procedures to further improve their solution quality.

**Table 1:** Related IDRT studies with solution methods.

| | P/F[1] | Obj.[2] | ST[3] | MD[4] | HV[5] | MR[6] | Exact/ heuristic | Solution method | Test instances |
|---|---|---|---|---|---|---|---|---|---|
| Liaw et al. (1996) | P | a | ✓ | | | | Heuristic | Online and offline decision support systems | Up to 120 requests |
| Hickman and Blume (2001) | P | d | ✓ | | | ✓ | Heuristic | Two-stage heuristic | Multiple transit lines, 3588 one-way requests |
| Aldaihani and Dessouky (2003) | P | a or b | ✓ | | | | Heuristic | Tabu search | Multiple transit lines, up to 155 requests |
| Häll et al. (2009) | P | a | | | | ✓ | Exact | MILP solver | One transit line, 4 requests |
| Ghilas et al. (2016) | F | ab | ✓ | ✓ | ✓ | | Metaheuristic | ALNS | One fixed line, up to 100 requests |
| Posada et al. (2017) | P | ab | ✓ | | ✓ | ✓ | Exact | MILP solver | One transit line, up to 6 requests |
| Ghilas et al. (2018) | F | ab | ✓ | ✓ | ✓ | | Exact | Branch-and-price algorithm | Two transit lines, up to 50 requests |
| Posada and Häll (2020) | P | ab | ✓ | | ✓ | ✓ | Metaheuristic | ALNS | One transit line, up to 145 requests |
| Molenbruch et al. (2021) | P | ab | ✓ | ✓ | ✓ | ✓ | Metaheuristic | LNS | Multiple transit lines, 200 requests |
| Melis et al. (2024) | P | c | | | | | Heuristic | Insertion-based heuristic | Multiple transit lines, up to 2000 requests |
| Fang et al. (2024) | P | ac | ✓ | ✓ | ✓ | ✓ | Exact | MILP solver | Two transit lines with up to 20 requests |
| This study | P | ac | ✓ | ✓ | ✓ | ✓ | Metaheuristic | Hybrid LNS | Multiple transit lines, up to 100 requests |

Notes: 1. P: passenger transport; F: freight transport. 2. Objective function: a. travel time (costs) of demand-responsive vehicles; b. travel time (costs) of transit line; c. customers' journey/travel time; d. generalized cost of customer inconvenience. 3. Synchronized transfer. 4. Multiple depots. 5. Heterogeneous vehicles. 6. Maximum ride (travel) time constraint. **Note that all the other studies do not consider electric vehicles.**

Table 1 summarizes the related IDRT studies discussed in this section and compares them with this study. Among the studies reviewed, this study is the only one to consider a weighted sum of both operator and customer costs and includes electric vehicles in the IDRT services. Furthermore, we extend IDRT in the context of Electric Vehicles (EVs), necessitating the development of efficient algorithms to jointly optimize vehicle routing and charging schedules that involve complicated trade-offs between operator costs and customer inconvenience, while satisfying multiple interdependent constraints related to vehicle routing and charging operations. Existing studies are limited on the EDARP without integrating fixed-route transit services. Incorporating a partial-recharge policy under capacitated charging-station constraints significantly complicates the problem, as the timing, location, and duration of each recharge event become interdependent across all buses. For a more in-depth algorithmic discussion, we refer readers to (Lam et al., 2022; Ma et al., 2024).

Regarding overall solution algorithms, three studies applied (A)LNS to their developed IDRT services, utilizing tailored destroy and repair operators for the specific characteristics of transit integration. As Turkeš et al. (2021) reported only a 0.14% improvement in solution quality when incorporating the adaptive layer of ALNS, we focus on developing an LNS-based algorithm. To improve the algorithm's performance, a hybridization approach is adopted together with the DA. We incorporate several novel local search operators, tailored for addressing the compact EIDARP. These local search operators enable a significant improvement in solution quality by balancing customer travel times and bus travel times.

## 3   The electric integrated dial-a-ride problem

In this section, we present the problem description of the EIDARP and summarize the characteristics of the EIDARP model. The detailed MILP formulation of the EIDARP, along with its notations, is presented in Appendices A and B.

### 3.1   Problem description

Consider a demand-responsive bus operator providing an integrated dial-a-ride service using a fleet of electric vehicles (buses). This service combines three travel modes: walking, demand-responsive bus (hereafter referred to as bus), and fixed-route transit services (referred to as transit hereafter). The bus fleet is heterogeneous in terms of passenger seats and battery size. Each bus, if used, needs to depart and return to its designated depot. The state of charge (SoC) of buses must be maintained within a certain range during the operating hours. Charging is restricted to operator-owned charging stations, and buses can only recharge when there are no passengers on board. We assume a constant discharging rate of buses when traveling and a constant charging rate when recharging. Charging operations are subject to capacity constraints at the charging stations. This is modeled at each individual charger, meaning that charging operations cannot overlap at any charger. Transit services are operated according to scheduled timetables with a constant dwelling time at each station. We assume the capacity of trains to be sufficiently large compared with the number of ride requests. Customers need to reserve the service in advance, specifying their origin, destination, and desired pickup or drop-off time windows. The operator informs customers of the details regarding the involved transport modes, including bus departure and arrival times, as well as transit (train) departure and arrival times, if applicable. Service options (customers' trips) from customers' origins to destinations can be one of the following five options:
1) Walk + transit + walk
2) Bus + transit + walk
3) Walk + transit + bus
4) Bus + transit + bus
5) Bus only

Note that transfers from one bus to another are not allowed. Customers are picked up and dropped off directly at their origin and/or destination if they are served by buses. For transfers between buses and transit stations, we assume that bus stops are located near their corresponding transit stations and that the additional walking time between them can be neglected. For customers whose origin or destination overlaps with transit stations, their walking time is neglected. The first four options distinguish EIDARP from DARP, where customers' first-mile/last-mile connections to transit stations are served by bus or walking. When buses serve customers' first-/last-mile, they are required to arrive at transit stations within a maximum waiting time (e.g., 10 minutes) before transit departures (first mile) or after transit arrivals (last mile) to ensure reasonable intermodal waiting times. For customer convenience, each customer is guaranteed a maximum travel time (door-to-door) characterized by a detour factor, a maximum walking distance (if walking is involved), and a maximum waiting time at transit stations (if using transit services). Customers may be rejected and incur a high penalty cost for the operator. Given a set of requests, the bus operator's objective is to minimize the weighted sum of bus and customer travel times, as well as the penalties for unserved requests.

### 3.2 Modeling EIDARP

The EIDARP is modeled on a directed graph $\mathcal{G} = (N, \mathcal{A})$, in which $N$ is the set of nodes and $\mathcal{A}$ is the set of arcs. $N$ is decomposed into five subsets related to different node types: customer origins, customer destinations, transit stations (stops), chargers, and depots. Let $R$ be the set of customer requests, i.e., $R = \{1,\ldots,n\}$. Customers' origin nodes and destination nodes are denoted as $P = \{1,\ldots,n\}$, and $D = \{n+1,\ldots,2n\}$, respectively. The transit network is modeled as a **departure-expanded network** (similar to a time-expanded graph, see Appendix C for a more detailed description), in which each transit station $f$ is duplicated according to its scheduled transit departures. In this structure, each duplicated node is associated with specific departure and arrival times. Let $\mathcal{F}_l$ and $\mathcal{D}_l$ be the set of physical transit stations and the set of departures for a transit line $l \in \mathcal{L}$. The set of transit (dummy) nodes is defined as $G = \{g_d^f | f \in \mathcal{F}_l, d \in \mathcal{D}_l, l \in \mathcal{L}\}$, where $g_d^f$ represents a transit node of physical station $f$ on transit line $l$ at departure $d$. The set of charger nodes consists of a set of charger dummy nodes $S$, where each charger (with a capacity of one) is duplicated with $m$ ordered dummies to allow for a maximum of $m$ visits from each bus. Note that we model the charging station capacity constraint at each individual charger based on the node duplication approach (Ma et al., 2024). Each charger node $s \in S$ has a specific charging speed (power) $\alpha_s$. The set of origin and destination depots is denoted as $O$ and $\bar{O}$, respectively.

The directed graph $\mathcal{G}$ comprises three subgraphs to model passenger/vehicle flows on bus, walk, and transit networks. Bus subgraph $\mathcal{G}_B$ is defined as $\mathcal{G}_B = (P \cup D \cup G \cup S \cup O \cup \bar{O}, \mathcal{A}_B)$, walking subgraph $\mathcal{G}_R = (P \cup D \cup G, \mathcal{A}_R)$, and departure-expanded transit subgraph $\mathcal{G}_T = (G, \mathcal{A}_G)$, where $\mathcal{A}_B$ is the set of bus arcs containing directed arcs between any pair of nodes in $\mathcal{G}_B$, and $\mathcal{A}_R$ is the set of walking arcs connecting customers' origins/destinations from/to transit nodes within customers' maximum walking distance. For the departure-expanded transit subgraph $\mathcal{G}_T$, the set of transit pairs (TPs) $\mathcal{A}_G$ represents connections between pairs of transit nodes in $G$ if a feasible path exists within a maximum transfer time $\gamma$ (e.g., 10 minutes). The algorithm to generate the TPs is provided in Appendix C. Each arc $(i,j) \in \mathcal{A}_G$, referred to as an origin-destination TP, is associated with a departure time $\underline{\theta}_i$ and arrival time $\bar{\theta}_j$. Each arc is assigned a distance $c_{ij}$. An illustrative example of the departure-expanded transit graph is presented in Appendix C. Arc travel times on different types of subgraphs are based on the speed of their transportation mode. We denote $t_{ij}$, $t_{ij}^w$, $t_{ij}^g$ as bus travel time, customer walking time, and transit travel time from node $i$ to $j$, respectively. In the subgraph $\mathcal{G}_B$, each node $i$ is associated with a service time $\mu_i$ and a time window $[e_i, l_i]$. The time window associated with transit node $i \in G$ depends on the user-defined maximum transfer time at transit stations, described in Section 4.1.

The operator's bus fleet is denoted by the set $K$, where each bus $k \in K$ is characterized by its passenger capacity $Q^k$ and energy consumption rate $\beta^k$. Each bus $k$ is assigned to an origin depot $o_k \in O$ and a destination depot $\bar{o}_k \in \bar{O}_k$. As buses are equipped with different battery capacities, they have dedicated minimum and maximum state of charge levels during operation, denoted as $E_{min}^k$ and $E_{max}^k$, respectively. The objective function for the EIDARP is presented in Eq. (1).

$$\text{Min } \lambda_1 \sum_{(i,j) \in A_B} \sum_{k \in K} t_{ij} x_{ij}^k + \lambda_2 \sum_{r \in R} \bar{L}_r + \omega \sum_{r \in R} (1 - v_r) \quad (1)$$

The first term of Eq. (1) represents the total bus travel times, where $x_{ij}^k$ is the decision variable indicating if an arc $(i,j)$ is traversed by bus $k$. The second term represents the total customer in-vehicle travel time and walking time, where $\bar{L}_r$ is calculated as the summation of customer in-bus time, in-transit (train) time, and walking time. Customer waiting times are not directly minimized but controlled via constraints. Note that we impose two constraints to limit customers' waiting times at transit stations and their maximum travel time. First, we limit the maximum transfer (waiting) time (user-defined), $\gamma$, when traveling between bus and transit vehicles. Second, we limit the maximum travel time of customers, characterized as their direct bus travel time multiplied by a detour factor. Thereby the computational efficiency of solving the MILP problem exactly can be significantly enhanced, compared to directly incorporating customers' door-to-door travel times into the objective function. Note that an alternative formulation and solution algorithm, which considers the weighted sum of total bus travel time and total excess user travel time, can be adopted (Su et al., 2023).

These two terms are weighted by $\lambda_1$ and $\lambda_2$, respectively. The last term penalizes customer rejection, where $\omega$ is the penalty associated with one customer rejection and $v_r$ is the binary decision variable determining whether customer $r$ is rejected.

The main characteristics of the considered service or operational constraints for the EIDARP are listed below.

(1) Each bus $k$ starts at its origin depot $o_k$ and ends at its destination depot $\bar{o}_k$.
(2) The number of customers on board for bus $k$ should not exceed its capacity $Q^k$.
(3) The SoC of a bus $k$ should be within the range of $E_{min}^k$ and $E_{max}^k$.
(4) Charging operations can be realized when no customer is on board.
(5) Each charger serves only one vehicle at a time (capacitated charging station constraints).
(6) Bus service times at a visited node need to be within the defined time window $[e_i, l_i]$.
(7) Walking is allowed when the walking time between customers' origins/destinations and transit stations is within the maximum walking time $t_{max}^w$.
(8) Customer's travel time from their origin to destination, denoted as $L_r$, computed as the difference between their arrival time and their departure time, cannot exceed their maximum travel time $L_r^{max} = t_{r,n+r} \times \varphi$, where $\varphi$ is the user-defined detour factor and $t_{r,n+r}$ is the direct travel time from the customer origin to the destination by a bus. It is possible to apply a customer-specific detour factor, $\varphi_r$, to meet the heterogeneous needs of customers.
(9) When a customer transfers between a bus and a transit, the transfer (waiting) time at the transit station is no more than the defined maximum transfer (waiting) time $\gamma$.

The reader is referred to Appendix B for the mathematical formulations of the constraints.

## 4    Solution algorithm

In this section, we propose a hybrid LNS to solve the EIDARP. The hybrid LNS integrates tailored local search operators and DA acceptance criteria to balance between diversification and intensification. Note that in our algorithm, the term "hybridization" refers specifically to the integration of tailored local search operators into the LNS metaheuristic to solve the EIDARP. It does not refer to the hybridization between LNS and DA. We first describe the preprocessing procedure to set up the time windows and tighten the solution space in Section 4.1. The graph after the preprocessing serves as input for the hybrid LNS in Section 4.2. In this algorithm, only feasible solutions determined by the nine-step solution evaluation scheme are accepted when a new solution is constructed. The solution feasibility is checked by the proposed nine-step evaluation scheme (Section 4.5). The feasibility evaluation scheme is particularly challenging due to the incorporation of multiple travel options for customers' journeys. The five service options defined in Section 2 introduce significant differences from DARP, impacting both the solution evaluation process and the design of destroy, repair, and local search (Sections 4.3-4.4) operators. Finally, we develop an efficient charging scheduling algorithm to schedule charging operations for the EIDARP (Section 4.6).

### 4.1 Preprocessing and representation of customer journeys

As introduced in Section 3.1, each transit node in the departure-expanded transit graph $\mathcal{G}_T$ is associated with departure and arrival times of a transit service. Given the predefined maximum allowed transfer time, we can establish the service time windows at the transit nodes. These time windows enable us to determine, for each request, a set of feasible TPs without violating the time windows and the maximum ride time of customers. These TPs serve as potential entry-exit transit station candidates for customers who can use transit as a part of their journeys. We denote $t_{ri}^w$ and $t_{ri}$ the waking and bus travel times from customer origin $r$ to the transit node $i$, respectively. $t_{j,n+r}^w$ and $t_{j,n+r}$ are the walking and bus travel times from transit node $j$ to customer destination $n+r$, respectively. $\underline{\theta}_i$ and $\bar{\theta}_j$ denote (transit/train) departure time at transit node $i$ (first-mile connection) and arrival time at transit node $j$, respectively. $\gamma$ is the maximum waiting/transfer time. Let $\mathcal{T}_r$ denote the set of TPs of customer $r \in R$. A TP $(i,j)$ can be added to $\mathcal{T}_r$ if it satisfies one of the following criteria given customer $r'$s origin and destination time windows $[e_r, l_r]$ and $[e_{n+r}, l_{n+r}]$, respectively.

- at origin: $e_r + t_{ri}^w \leq \underline{\theta}_i$, $l_r + t_{ri}^w \geq \underline{\theta}_i - \gamma$, and $t_{ri}^w$ is within the maximum walking time; at destination: $\bar{\theta}_j + t_{j,n+r}^w \leq l_{n+r}$, $\bar{\theta}_j + t_{j,n+r}^w + \gamma \geq e_{n+r}$ and $t_{j,n+r}^w$ is within the maximum walking time (travel option 1)
- at origin: $e_r + t_{ri} \leq \underline{\theta}_i$ and $l_r + t_{ri} \geq \underline{\theta}_i - \gamma$; at destination: $\bar{\theta}_j + t_{j,n+r}^w \leq l_{n+r}$, $\bar{\theta}_j + t_{j,n+r}^w + \gamma \geq e_{n+r}$ and $t_{j,n+r}^w$ is within the maximum walking time (travel option 2)
- at origin: $e_r + t_{ri}^w \leq \underline{\theta}_i$, $l_r + t_{ri}^w \geq \underline{\theta}_i - \gamma$, and $t_{ri}^w$ is within the maximum walking time; at destination: $\bar{\theta}_j + t_{j,n+r} \leq l_{n+r}$ and $\bar{\theta}_j + t_{j,n+r} + \gamma \geq e_{n+r}$ (travel option 3)
- at origin: $e_r + t_{ri} \leq \underline{\theta}_i$ and $l_r + t_{ri} \geq \underline{\theta}_i - \gamma$; and at destination: $\bar{\theta}_j + t_{j,n+r} \leq l_{n+r}$ and $\bar{\theta}_j + t_{j,n+r} + \gamma \geq e_{n+r}$ (travel option 4)

Let $A_i^k$ and $B_i^k$ denote the arrival time and starting time of service of a bus $k$ at transit node $i$, respectively. When arriving at a transit node, a bus's arrival time or starting time of service needs to satisfy the following time window constraints:

1) If a bus $k$ has customers transferring to transit services, its arrival time needs to satisfy $\underline{\theta}_i - \gamma \leq A_i^k \leq \underline{\theta}_i$.
2) If a transit node $i$ has customers transferring from transit services to bus services on bus $k$, the starting time of service needs to satisfy $\bar{\theta}_i \leq B_i^k \leq \bar{\theta}_i + \gamma$.
3) If transit node $i$ has both customers transferring between bus and transit service on bus $k$, it needs to satisfy $\underline{\theta}_i - \gamma \leq A_i^k \leq \underline{\theta}_i$ and $\bar{\theta}_i \leq B_i^k \leq \bar{\theta}_i + \gamma$.

We initialize the time windows at transit node $i$ in the preprocessing step as $[\underline{\theta}_i - \gamma, \bar{\theta}_i + \gamma]$. The time windows of customers' origins or destinations are input data provided by customers, whereas the time windows for depots and chargers are defined as $[0, t_{end}]$, where $t_{end}$ is the end of the service operation horizon. We apply the time-window tightening and arc elimination procedures for bus nodes using the method proposed by Fang et al. (2024).

For the DARP, a customer's journey is simply a sequence of nodes visited by a bus connecting their origin and destination. In (E)IDARP, however, a customer's journey might include a transit leg requiring additional bus service for first- and/or last-mile bus service. Figure 1 illustrates a solution where customers utilize different service options. Let $i^+$ and $i^-$ denote the pickup and drop-off nodes of customer $i$, respectively. The customer journey is stored as a sequence of visited nodes, and the visited transit nodes are associated with labels indicating whether there are first-mile or last-mile connections and the connected transport mode (walking or bus) used to reach them. In the algorithm, we store the travel times of each trip leg (first mile, transit leg, last mile) of customers' journeys in a matrix, which allows us to extract this information efficiently for the solution feasibility evaluation when updating part of the customer journeys. In this example, two bus routes are used. Bus 1 leaves from depot $o^1$ and brings customers 1 and 3 to transit station $g_1$, while bus 2 starts from $o^2$ to serve customers 2 and 3's last miles from transit station $g_2$. The representation of each customer's journey is as follows.

1) Customer 1 (bus-transit-walk): $1^+ - 3^+ - 4^+ - 4^- - g_1 - g_2 - 1^-$
2) Customer 2 (walk-transit-bus): $2^+ - g_1 - g_2 - 2^-$
3) Customer 3 (bus-transit-bus): $3^+ - 4^+ - 4^- - g_1 - g_2 - 2^- - 3^-$
4) Customer 4 (bus only): $4^+ - 4^-$

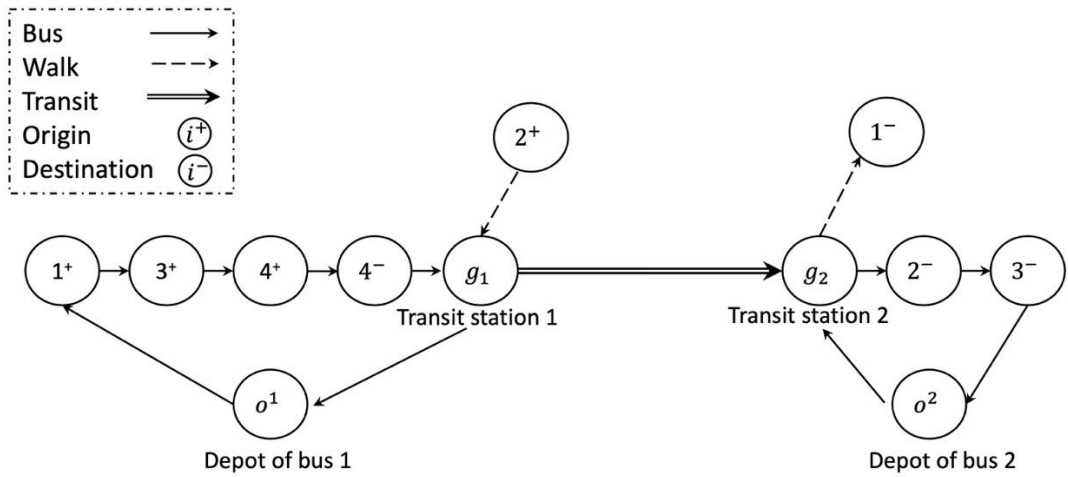

**Figure 1:** Representation of customer journeys and bus routes for IDARP services involving different service options.

## 4.2 The hybrid large neighborhood search

The framework of the hybrid LNS is presented in Algorithm 1. First, an initial solution $s_{init}$ is constructed by prioritizing customers who can use transit services (i.e., those with a TP set that is not empty). In doing so, bus travel costs can be reduced as part of customers' journeys are served by the transit service. Customers who are eligible to use transit services are sorted by their earliest time window at their origin node. They are then inserted one by one, with travel options 1 through 4 being checked in ascending order until a feasible option is found. If none of these four options is feasible, the customer is assigned to the bus-only service (travel option 5). New bus routes are generated as needed.

Once the $s_{init}$ is constructed, line 2 sets both the incumbent solution $s$ and the current best solution $s_{best}$ as $s_{init}$. In the next step, a destroy operator (i.e., removing customers from current bus routes) and a repair operator (i.e., inserting unserved customers into current bus routes) are randomly selected and applied to $s$, yielding $s'$. For the destroy operators, $\xi_{max} \in (0,1)$ is a parameter that sets the degree of destruction of a current solution and determines the size of the search neighborhood (line 6). When $\xi_{max}$ is small (e.g., 0.05), the algorithm cannot move to a farther search area to escape from the local optimum. However, when $\xi_{max}$ is large (e.g., 0.5), it may destroy a good solution and require more computational time to find others (Ropke and Pisinger, 2006). The value of $\xi_{max}$ needs to be tuned to ensure good performance. We tune its value and other parameters of the algorithm in Section 5.2. The number of removed customers $n^{destroy}$ is randomly selected between 1 and $\lceil n \times \xi_{max} \rceil$ where $n$ is the number of customers. Subsequently, a sequence of local search operators is applied to $s'$ to obtain a local optimal solution $s''$ (line 10). Local search procedure is only activated when $s'$ is promising, i.e., $c(s') \leq \alpha c(s_{best})$, where a smaller $\alpha$ ($\alpha > 1$) suggests less-frequent local search to save computational times.

A DA acceptance criterion decides whether $s''$ should be accepted as the next incumbent solution. At each iteration, $s_{best}$ is updated if the objective function value $c(s'') \leq c(s_{best})$ (lines 16-18). The acceptance criterion is based on the DA approach. $s''$ is accepted at the next iteration if $c(s'') \leq c(s_{best}) + T$, where the temperature $T$ is updated at every iteration as $T = T - T_{max}/T_{red}$. $T_{red}$ is a cooling parameter. The maximum temperature $T_{max}$ depends on the average distance $\bar{c}$ of the arcs of the graph, i.e., $T_{max} = t_{max} \times \bar{c}$ (Braekers et al., 2014), where $t_{max}$ is a parameter to be set (see Section 5.2). $T_{red}$ is the cooling temperature that decides the cooling speed of $T$. $T$ is reset to $T_{max}$ when $T \leq 0$ (line 12). Note that we adopt the standard DA approach in which the temperature parameter $T$ is simply restarted to $T_{max}$ when $T$ becomes negative (Santini et al., 2018). An additional restart mechanism can be incorporated into the DA algorithm when the algorithm stagnates at local optima after a certain number of iterations (Braekers et al., 2014). The algorithm terminates when the maximum iteration is reached, and $s_{best}$ is returned as the final solution.

| Algorithm 1. Hybrid LNS. |
| --- |
| Input: $\mathcal{G}$, $K$, time windows, $\lambda_1, \lambda_2, \omega$, and $\varphi$ |
| Output: $s_{best}$ |
| 1: Construct $s_{init}$ |
| 2: $s, s_{best} \leftarrow s_{init}$; $i = 1$ |
| 3: **repeat** |
| 4:     randomly select a destroy operator $d_i \in \{d_{random}, d_{worst}, d_{relate}, d_{route}\}$ |
| 5:     randomly select a repair operator $r_i \in \{r_{random}, r_{greedy}, r_{regret}, r_{Tpriority}, r_{TP}\}$ |
| 6:     $n^{destroy} = rand(1, \lceil n \times \xi_{max} \rceil)$ |
| 7:     ($s$, unserved customers) $\leftarrow d_i(s, n^{destroy})$ |
| 8:     $s' \leftarrow r_i(s,$ unserved customers$)$ |
| 9:     **if** $c(s') < \alpha \times c(s_{best})$ |
| 10:         $s'' \leftarrow LocalSearch(s')$ |
| 11:     **end if** |
| 12:     $T = T - T_{max}/T_{red}$; if $T \leq 0$, set $T = T_{max}$ |
| 13:     **if** $c(s'') < c(s_{best}) + T$ |
| 14:         $s \leftarrow s''$ |
| 15:     **end if** |
| 16:     **if** $c(s'') < c(s_{best})$ |

17:         $s_{best} \leftarrow s''$
18:     **end if**
19:     $i := i + 1$
20: **until** $i = n_{iter}$
21: **return** $s_{best}$

## 4.3 Destroy and repair operators

The shaking step consists of destroy and repair procedures (operators), with most operators having been extensively applied in various vehicle routing problems and their variants (Mara et al., 2022). We propose three tailored destroy/repair operators, $d_{relate}, r_{Tpriority}$, and $r_{TP}$, for the EIDARP. When removing a customer, the involved bus routes are updated, and related walking legs to transit stations are removed. The set of proposed destroy operators includes:

- *Random removal* – $d_{random}$: randomly remove $n^{destroy}$ customers from the served customers.
- *Worst removal* – $d_{worst}$: remove the $n^{destroy}$ worst customers in terms of objective function value.
- *Route removal* – $d_{route}$: randomly remove a bus route. If the route serves a customer only for their first- or last-mile segment (travel options 2, 3, or 4), that customer's entire trip is removed from both used transit and bus services.
- *Related removal* – $d_{related}$: remove $n^{remove}$ most related customers from the served customer list in solution $s$. Eq. (2) defines the relatedness of customers $i$ and $j$, where $c_{max}$ represents the distance of the longest arc in $\mathcal{G}$, and $t_{max}$ is the end of the planning horizon. The first term in Eq. (2) represents the spatial relatedness of their origins and destinations, and the second term computes the temporal relatedness of their origins and destinations. For customers' journeys involving transit legs, the relatedness of customer $i$'s and $j$'s TPs are incorporated as a third term by calculating the ratio of shared TPs, i.e. $|\mathcal{T}_i \cap \mathcal{T}_j|$ divided by the total number of TPs, i.e. $|\cup_{r \in R} \mathcal{T}_r|$. Note that $abs(A_i - A_j)$ represented the absolute value of $A_i - A_j$, while $|\mathcal{T}_i \cap \mathcal{T}_j|$ represents the number of elements in the set $\mathcal{T}_i \cap \mathcal{T}_j$.

$$Relatedness_{ij} = \left( \frac{c_{ij} + c_{n+i,n+j}}{2c_{max}} + \frac{abs(B_i - B_j) + abs(A_i - A_j)}{2t_{max}} + \frac{|\mathcal{T}_i \cap \mathcal{T}_j|}{|\cup_{r \in R} \mathcal{T}_r|} \right)^{-1} \quad (2)$$

Repairing a customer in the algorithm aims to find the best option and the resulting route insertion positions with the least cost (i.e., the best insertion position). In DARP, this process only involves iterating over all routes to determine the best insertion positions (travel option 5). However, for EIDARP, the complexity is significantly higher due to the need to evaluate all five different travel options for each customer. In the case of EIDARP, searching for the best position involves finding the best travel option along with its associated best insertion positions. For customers using transit service in travel options 1–4, we consider their first-mile legs first, followed by their last-mile legs, as the customers' maximum travel time constraint is evaluated and adjusted at their destination (detailed in Section 4.5).

- *Random repair* – $r_{random}$: repair unserved customers in a random order at the lowest-cost insertion positions, i.e., positions that result in the smallest increase in the objective function value.
- *Greedy repair* – $r_{greedy}$: iteratively select and insert a customer whose best insertion position has the lowest cost among the unserved customers.
- *Regret repair* – $r_{k-regret}$: iteratively select and insert a customer who has the highest regret value, which is calculated as the cumulative difference between the objective function value at the best insert position and the $k^{th}$-best insertion position. In this algorithm, $k$ is chosen as 2 or 3 randomly (Ropke and Pisinger, 2006).
- *Transit priority repair – $r_{Tpriority}$*: repairs unserved customers by those who have the option of using trains (i.e., those who have a non-empty set of TPs), thereby reducing the bus travel costs in the objective function.
- *TP priority repair – $r_{TP}$*: repairs multiple customers at the same time if they share the same TP(s) as they may have a higher likelihood to share the same first-mile/last-mile bus trips, leading to a reduction in objective function value.

### 4.4 Local search operators

Local search operators (e.g., *2-opt, 4-opt, relocate*) are widely used in DARP, PDP, and vehicle routing problems. These search operators aim to reduce bus travel times (costs) from a bus operator's perspective. However, our EIDARP aims to minimize both customer and bus travel times as a weighted sum. Our proposed local search operators focus on the trade-offs between the two costs. Five local search operators are applied as described below:
- *Bus exchange*: exchange a bus with another bus starting from a different depot, provided that such an exchange results in a lower objective function value while maintaining the sequences of visited customers unchanged.
- *Exchange TPs of customers*: for customers served by transit service, replace each customer's TP with another available TP candidate in their respective TP set. The first-mile and last-mile bus routes or walking legs might also be changed if a better solution is found. Update the solutions if the objective function value is reduced.
- *Replace TP by bus*: for customers served by transit, evaluate each customer individually to determine whether switching them to a door-to-door bus service reduces the objective function value. If so, retain the newly obtained solution.
- *Re-assign first-/last-mile bus service*: for customers served by transit, explore alternative bus routes for first- or last-mile connections. Retain the new route if the objective function value is improved.
- *Replace walk with bus service*: for customers' first-/last-mile legs served by walking, replace their walking arcs with a bus. Retain this change if it reduces the objective function value.

### 4.5 Feasibility evaluation for the integrated bus route

To address the specific requirements of IDARP, we extend the eight-step evaluation scheme of Parragh et al. (2010) to a nine-step one to check the feasibility of bus routes and the involved customers' journeys. The evaluation scheme ignores the energy feasibility of a bus route. Once the considered bus route is deemed feasible (regardless of energy feasibility), we verify whether it is also feasible in terms of energy. If not, we execute the energy-repair procedure described in Section 4.6 to evaluate the overall feasibility. The nine-step feasibility evaluation is used to evaluate the feasibility of bus routes. Similar to the classical DARP's eight-step feasibility evaluation scheme, vehicle capacity and time window constraints must be satisfied for a bus in the IDARP service. However, for the customer's maximum travel time (or maximum ride time in classical DARP) constraint, the feasibility evaluation and repair process differ. Each bus route in an IDARP might have customers who need to be dropped off at transit stations, then continue their journeys to their destination via transit + bus or transit + walk, or who might have completed the first part of their journey by another bus and transit, and continue their journeys on this bus to their destination (last-mile service). Consequently, when evaluating the IDARP customers on a bus route, the travel time of customers using travel options involving transit (options 2-4: bus-transit-walk, bus-transit-bus, and walk-transit-bus) needs to take into account the travel times on the subsequent or preceding trip legs connected to their origin/destination to obtain customer travel time from origin to destination. Note that if a bus route has no customers using the bus-transit-bus option, it is evaluated independently using the nine-step evaluation scheme. Otherwise, the feasibility evaluation needs to verify the maximum travel time feasibility for customers using the bus-transit-bus option, which involves other corresponding bus routes. The proposed evaluation scheme aims to address this more complex routing problem, ensuring that bus capacity, time windows, and customer travel time constraints are satisfied.

Given a bus route $\sigma$ starting and terminating at its designated depot, we define the following indicators ($\sigma$ is dropped) to be used in the nine-step evaluation scheme:

- $A_i$: Arrival time at node $i$, $A_i = D_{i-1} + t_{i-1,i}$
- $B_i$: Beginning time of service at node $i$, $B_i = \max(A_i, e_i)$
- $W_i$: Waiting time at node $i$, $W_i = B_i - A_i$
- $D_i$: Departure time at node $i$, $D_i = B_i + \mu_i$
- $q_i$: Passenger load at node $i$
- $L_i$: Customer $i$'s travel time between origin and destination. As mentioned above, for customers using service options 2-4, we need to retrieve their travel times on their subsequent or preceding trip legs to obtain the customer's travel times between their origin and destination.
- $F_i$: Forward slack time at node $i$ (described later)
- $\rho_i$: Waiting time violation at a transit node $i \in G_\sigma$ if customers transferring from the bus to a transit node $i$, $\rho_i := \left(\underline{\theta_i} - \gamma - A_i\right)^{+}$[1], where $G_\sigma$ is the set of transit nodes on the route for the first-mile transfers.

Note that in our EIDARP, we allow earlier arrivals than $e_i$ (i.e., $A_i < e_i$) at any node, with a waiting time of $W_i$, except at transit nodes where customers transfer from the bus. This is because bus arrival times at transit nodes must meet the customer's maximum transfer time constraints (i.e., $\underline{\theta_i} - \gamma \leq A_i \leq \underline{\theta_i}$). We introduce $\rho_i$ to calculate waiting time violations at transit nodes for the first mile service during the feasibility check in the nine-step evaluation scheme (described later). Note that for the last-mile services (i.e., bus-transit-bus or walk-transit-bus), buses are allowed to arrive earlier than $e_i$ at transit stations, as this will not affect customers' maximum waiting times at transit stations. Customers' arrival times at transit stations for their last-mile transfers are based on their trains' arrival times. We drop the index of bus $k$ for the above-mentioned variables in the following description for the feasibility evaluation and vehicle recharge scheduling.

For the IDARP, the forward slack times, $F_i$, are calculated as shown in Eq. (3), which accommodates the different travel options used by customers on the bus route.

$$F_i = \min_{i \leq j \leq \bar{q}} \left\{ \sum_{i < p \leq j} W_p + \varsigma_j \right\} \quad (3)$$

where $\varsigma_j$ equals:

a. $\left(\underline{\theta_j} - \max(A_j, \underline{\theta_j} - \gamma)\right)^{+}$, if $j$ is a transit node with customers transferring from bus to transit service.
b. $\left(\min\{l_j - B_j, L_r^{max} - L_j\}\right)^{+}$, if $j$ is a destination node of customer $r$.
c. $\left(l_j - B_j\right)^{+}$, otherwise.

$\bar{q}$ denotes the end node of the corresponding route.

---

[1] $(x)^{+} = \max\{0, x\}$

This equation is a generalized formulation for calculating $F_i$ in which classical DARP is a special case. The first term of Eq. (3) is the cumulative waiting time from $i$ to the end $\bar{q}$ of the route, identical to $F_i$ calculation for the DARP (Cordeau and Laporte, 2003). The second term, $\varsigma_j$, is calculated differently depending on the node being visited. If the node is a transit node with customers transferring from a bus to it (case a), the additional slack time equals $\left(\underline{\theta}_j - \max(A_j, \underline{\theta}_j - \gamma)\right)^+$, where $\underline{\theta}_j$ is the departure time of the train/transit to be connected by the bus. If the node is a destination node of customer $r$ (case b), it captures the maximum travel time constraint of the customer (whose destination node is $j$) and the latest time window at $j$, identical to the DARP. In case b, $L_j$ denotes the travel time for customer $r$. Our algorithm retrieves information of arrival or departure times on the other trip legs when transit or another bus is used as part of customer $r$'s journey in order to calculate $L_j$. Otherwise, the additional forward slack time is calculated as $\left(l_j - B_j\right)^+$ (case c). Eq. (3) considers transit departure synchronization and maximum waiting time constraints at transit stations.

Note that $\rho_i$ represents the waiting time violation when arriving at a transfer node $i$ for first-mile transfers. If customers are transferring from the bus route to take a transit service, the bus's arrival time cannot be earlier than $\underline{\theta}_j - \gamma$. $\rho_i$ serves as an indicator to determine whether it is necessary to apply the travel time repair procedure (Step 7) in the nine-step evaluation scheme.

Algorithm 2 represents the nine-step evaluation scheme. If a bus route contains customers who use the bus-transit-bus option and whose last mile has not yet been inserted, the maximum travel time constraints for these customers will be temporarily disregarded and will then be verified on the buses that serve their last miles. Therefore, when inserting a customer using the bus-transit-bus option, the first mile is inserted first. If a feasible position is found for the first mile, the customer's last mile is inserted on another bus route. Then, the maximum travel time constraint for this customer is checked on the corresponding last-mile route, as well as the feasibility constraints for other customers. If the corresponding last-mile route is not feasible, we check other routes for the last-mile insertion of this customer. If no feasible insertion can be found for the last-mile insertion of this customer, then the customer is removed from both the first-mile and last-mile bus routes.

Consider a bus route $\sigma$ to be evaluated. The nine-step evaluation scheme begins with setting the departure time from the depot to the earliest time window (step 1). In Step 2, $A_i$, $B_i$, $W_i$, $D_i$, and $q_i$ are updated, and potential violations of time window, bus capacity, and transit departure synchronization are checked. Steps 3-6 delay departure time without violating maximum waiting time and maximum travel time. Note that the maximum travel time verification does not apply to customers for whom the last mile has not yet been inserted. If neither the maximum travel time nor the transit departure synchronization is violated, the process advances directly to Step 9; otherwise, Step 7 is applied to delay the departure time at every origin node and transit node. Step 7(a) computes the forward slack time $F_j$. Step 7(b) delays the departure time to minimize customer travel time. This step also ensures minimum waiting time at the transit station if the customer changes from bus to transit service at node $j$. Step 7(c) updates all the indicators after node $j$. For transit nodes $i' \in G_\sigma$, $\rho_{i'}$ is recalculated as $\rho_{i'} = \left(\rho_{i'} - F_j\right)^+$, reflecting the updated waiting time violations at node $i'$ after the departure from node $j$ is delayed by $F_j$. Additionally, the flag $f_\sigma^{delay}$ is set to $true$, indicating that the delay procedure has already been applied to this route $\sigma$. Step 7(d) checks whether their travel time at the destination node $i$, $L_i$, exceeds the maximum travel time of customers ($L_i^{max}$) and whether the maximum waiting times at transit stations for first-mile transfers are satisfied. If all $L_i \leq L_i^{max}$ and $\sum_{i' \in G_\sigma} \rho_{i'} = 0$, this route is feasible, then proceed to Step 9. Otherwise, Step 8 is applied to repair the travel time infeasibility for customers using the bus-transit-bus option. This step delays the departure time along the first-mile routes for customers whose last mile lies on the examined route, thereby reducing the travel times of the corresponding customers on the examined route. Consider a request $i$ on the examined route $\sigma$, whose first mile is served by route $\sigma'$. If $f_{\sigma'}^{delay} = false$, it means the departure delay procedure (Step 7) has not been applied to route $\sigma'$. Step 8 then utilizes the forward slack times to delay the departure time at the origin/transit nodes for route $\sigma'$ while ensuring all the constraints are satisfied. This will eventually reduce the travel times for request $i$ and other requests (on route $\sigma$) whose first miles are served by route $\sigma'$. After this travel time tightening procedure, if the maximum travel time constraints are satisfied for all the requests on the examined route, the route is feasible. Otherwise, it is infeasible. Note that computing the changes in violations of bus capacity, route duration, time windows, and travel time constraint is optional. These breaches are not penalized in our objective function. We keep them to facilitate comparison with the eight-step feasibility evaluation scheme of the DARP.

---

**Algorithm 2.** Nine-step evaluation scheme.

Input: A bus route $\sigma$ to be evaluated.
Output: feasible = true or false.

1: Set $D_0 := e_0, f_\sigma^{delay} := false$
2: Compute $A_i, W_i, B_i, D_i$, and $q_i$ for each node $i$ on the route. If $B_i > l_i$, or $q_i > Q$, or if the bus arrives later than the departure of to-be-connected transit at the transit node $i$, i.e., $A_i > \underline{\theta_i}$, GO TO STEP 9.
3: Compute $F_0$.
4: Set $D_0 := e_0 + \min\{F_0, \sum_{0<p<\bar{q}} W_p\}$.
5: Update $A_i, W_i, B_i, D_i$ for each node $i$ on the route, and calculate $\rho_{i'}$ for $i' \in G_\sigma$.
6: Compute $L_i$ for each request on the route, excluding the requests for which their last miles have not yet been inserted. If all $L_i \leq L_i^{max}$ and $\sum_{i' \in G_\sigma} \rho_{i'} = 0$ (no maximum waiting time constraint violations at the transit nodes $G_\sigma$ for bus-transit-bus or bus-transit-walk transfers), GO TO STEP 9.
7: For every node $j$ that is an origin node or a transit node with customers transferring from the bus to a transit service:
   a) Compute $F_j$.
   b) Set $W_j := W_j + F_j; B_j := A_j + W_j; D_j := B_j + \mu_j$.
   c) Update $A_i, W_i, B_i, D_i$ for each node $i$ that comes after $j$ in the route. For each transit node $i' \in \bar{G}_\sigma$ that comes after $j$ in the route, update $\rho_{i'}$ as $(\rho_{i'} - F_j)^+$.
   d) Update $L_i$ for each request $i$ whose destination is after $j$. Set $f_\sigma^{delay} := true$. If all $L_i \leq L_i^{max}$ for requests whose destinations lie after $j$ and $\sum_{i' \in G_\sigma} \rho_{i'} = 0$, GO TO STEP 9; otherwise, GO TO STEP 8.
8: For every node that is the destination node of a last-mile request $i$ using the bus-transit-bus option on the current route $\sigma$:
   Compute $L_i$ by retrieving the travel time information on its first mile. Let $\sigma'$ denote the bus route of the first mile of $i$. If $L_i > L_i^{max}$ and $f_{\sigma'}^{delay} = false$, then apply STEP 7 for the first-mile route $\sigma'$ to delay the departure time at each origin node/transit node. Set $f_{\sigma'}^{delay} := true$. If successful, update all $L_i$ for request $i$ and other corresponding requests whose first miles are served by route $\sigma'$.
   If $L_i \leq L_i^{max}$ for all the requests on route $\sigma$, route $\sigma$ is feasible; otherwise, it is infeasible.
9: Compute changes in violations of bus load, duration, time window, and customer travel time constraints.

---

As IDARP deals with multi-modal customer journeys, the differences between the nine-step scheme (IDARP) and the eight-step scheme (DARP) are summarized as follows.

a. IDARP extends Step 2 of the eight-step scheme by verifying if the bus's arrival time at a transit station exceeds the transit departure time (i.e., $A_i > \underline{\theta_i}$).
b. IDARP introduces $\rho_i$, the waiting time violation at transit station $i$ for customers' first-mile connections. In Step 6, $\rho_i$ is checked to determine whether a travel time repair procedure (Step 7) is necessary.
c. In Steps 3 and 7, the forward slack time $F_i$ calculation differs (Eq. (3)) from the DARP, which is a generalized formulation in the context of IDARP.
d. In Step 7, IDARP delays the route's departure time at both customer origins and transit stations when there are customers transferring from buses to transit on the route, while DARP can delay only the departure time at customer origins. In Step 7(b), IDARP updates $W_j$ as $W_j := W_j + F_j$ instead of $W_j := W_j + \min\{F_j, \sum_{j<p<\bar{q}} W_p\}$ at an origin node $j$ on the route. Cordeau and Laporte (2003) demonstrate that delaying the departure time at any node $i$ within $\sum_{i<p\leq\bar{q}} W_p$ does not impact the bus arrival times at any node, but further delays would increase the arrival times of the subsequent nodes after node $i$ as much. Let $\Delta_j$ be $F_j - \sum_{j<p<\bar{q}} W_p \geq 0$. As illustrated in Figure 2(a), we can delay the departure as much as $\min\{F_j, \sum_{j<p<\bar{q}} W_p\}$ at an origin node $j$ for the DARP to reduce the route duration while satisfying time window and travel time constraints. However, for IDARP, when customers' origins and destinations are served by different buses, using the forward slack times to delay the departure time at any node for one bus route will not affect the arrival time at any node of another route. This allows us to delay extra departure time $\Delta_j$ at an origin node $j$ by more than $\sum_{j<p\leq\bar{q}} W_p$ for the first-mile route, as it will not affect the arrival time at the destination $n+j$ on its corresponding last-mile route. In doing so, this extra departure time delay can still reduce customer $j$'s travel time (see Figure 2(b)).
e. Step 8 is relevant only for customers with the travel option 4 (bus + transit + bus), where the first-mile bus route has not executed the departure delay procedure (Step 7).

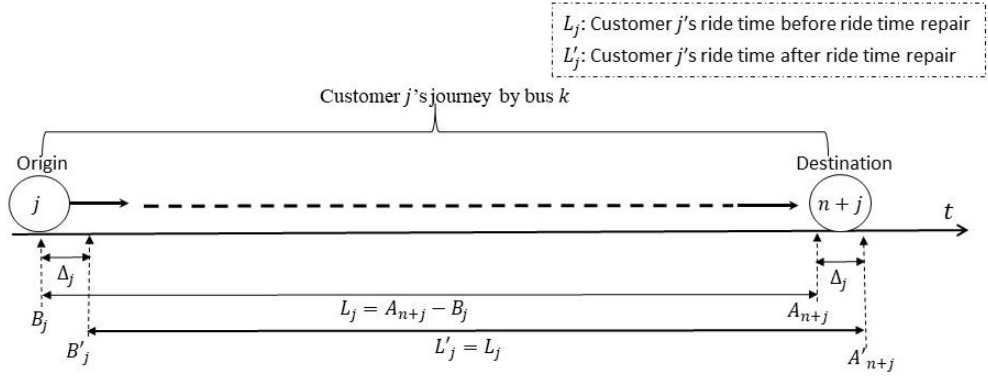

(a) DARP: Delay extra departure time $\Delta_j$ at an origin node $j$ would increase the arrival times at the destination node $n+j$ as much as $\Delta_j$.

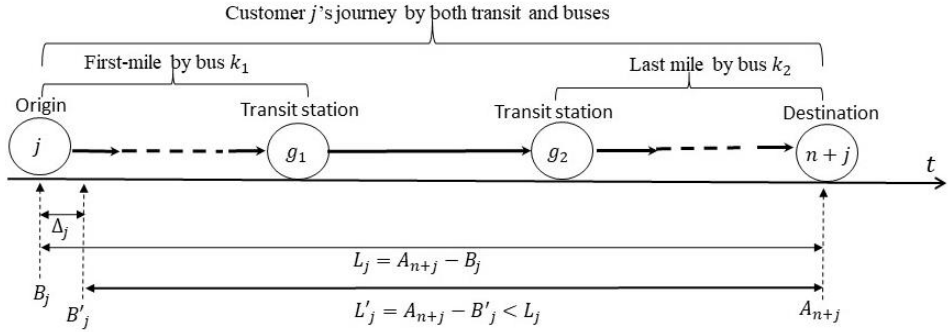

(b) IDARP (bus-transit-bus): Delaying extra departure time $\Delta_j$ at an origin node $j$ more than $\sum_{j<p\leq\bar{q}} W_p$ for the first mile can still reduce customer $j$'s travel time, hence repair customer $j$'s travel time infeasibility.

**Figure 2:** Further delay customer $j$'s departure time $\Delta_j$ at Step 7(b) for DARP and IDARP when the departure time at the node $j$ has already been delayed $\sum_{j<p<\bar{q}} W_p$.

### 4.6 Vehicle charging scheduling

The vehicle charging (recharge) scheduling procedure aims to insert recharging operations when a feasible route is energy-infeasible without recharging. This algorithm extends the method of Ma et al. (2024) to handle capacitated charging stations and to compute the available recharging time in the context of the EIDARP. Moreover, we adopt a similar charging scheduling strategy, as used in Bongiovanni et al. (2024), to schedule bus charging operations as early as possible and as much as possible from the first candidate position (the depot) to insert charging operations. Note that the proposed charging scheduling algorithm is a heuristic approach to efficiently find a good feasible charging schedule.

Let $\Delta_E$ denote the amount of energy violation of the current route, computed as $\Delta_E = E_{min} - (E_{init} - ec)$, where $E_{init}$ and $E_{min}$ are the initial and minimum SoC, respectively. $ec$ is the total energy consumption of the examined route. If $\Delta_E \leq 0$, the recharging scheduling is skipped. Otherwise, we execute this algorithm to repair its energy infeasibility. First, we identify the first position $i_{low}$ where the vehicle's SoC is below $E_{min}$, meaning that a recharge event must be inserted before $i_{low}$ with the amount $\Delta_{E_{low}} = E_{min} - E_{i_{low}}$, where $E_{i_{low}}$ is the SoC at position $i_{low}$. $\Delta_E$ and $\Delta_{E_{low}}$ represent the total required recharging amount and minimum required recharging amount, respectively. Each charging node $s \in S$ maintains a record of its charging events, including the IDs of visited buses, the start and end times of their charging operations, and their charging durations. Based on this information, the available (vacant) time intervals of $s$ can be calculated as the gaps between the recorded charging events. The recharging scheduling aims to find common time intervals that are feasible for both the bus route and the charging station.

The vehicle charging algorithm is described in Algorithm 3. It starts with seeking a feasible recharge between the first position and $i_{low}$. Let $I$ denote an ordered set of nodes of the route with no passengers onboard (i.e., $q_i = 0$). It starts by checking the nearest chargers (line 4). Line 5 computes the maximum potential recharging time $\delta_{v_is}$ at position $i$ if recharged at charging node $s$, which is determined by deducting $F_i$ from the travel time difference to $s$ and the service time at $s$ (i.e., $\mu_s$). $\Delta_E^s$ and $\Delta_{E_{low}}^s$ are updated by including extra energy consumption traveling to $s$, with the total (maximum) required charging time $\delta_{v_is}^{max}$ and the minimum required charging time $\delta_{v_is}^{min}$ calculated based on the charging power of $s$. If $\delta_{v_is} < \delta_{v_is}^{min}$, the position $i$ is infeasible and the algorithm moves to the next one (line 8); otherwise, the maximum available recharging time interval $\varepsilon_{route}^i$ at position $i$ is calculated (line 9). This interval is then compared with the available time intervals $\varepsilon_s$ at charger $s$. If a common time interval exists, but the duration $\delta$ cannot meet the minimum recharging time $\delta_{v_is}^{min}$, it goes to line 10 and searches for another one. If the duration $\delta$ is more than $\delta_{v_is}^{max}$, the recharging is considered successful with charging time $\delta_{v_is}^{max}$, and this charging event is added to $s$ (lines 13-14). If a common time interval $\delta$ is between $\delta_{v_is}^{min}$ and $\delta_{v_is}^{max}$, the recharging event is inserted with the charging time $\delta$, but the scheduling process continues from the next recharge position (Go To line 2) with updated $\Delta_E$, $\Delta_{E_{min}}$ and $i_{low}$ (lines 16-18). If $\Delta_E$ is still greater than 0 at the end of the algorithm, the recharge scheduling is considered unsuccessful.

---

**Algorithm 3. Recharging scheduling with capacitated charging stations.**

Input: $\sigma = \{v_0, \ldots, v_q\}, S, \Delta_E, i_{low}, \Delta_{E_{low}}$
Output: success = true or false

1:    $i_{init} \coloneqq 0$
2:    **for** $i$ in $I$
3:      Compute $F_i$ based on Eq. (3).
4:      **for** $s \in S$, where $S$ is sorted by the distance to $v_i$ in ascending order
5:        Compute $\delta_{v_is}$ as $\delta_{v_is} \coloneqq F_i - (t_{v_is} + t_{sv_{i+1}} - t_{v_iv_{i+1}}) - \mu_s$.
6:        Update the total required recharge amount and minimum required recharge amount as $\Delta_E^s \coloneqq \Delta_E + (c_{v_is} + c_{sv_{i+1}} - c_{v_iv_{i+1}}) \times \beta, \Delta_{E_{low}}^s \coloneqq \Delta_E + (c_{v_is} + c_{sv_{i+1}} - c_{v_iv_{i+1}}) \times \beta$.
7:        Compute the total required and minimum required charging time if recharged at $s$, $\delta_{v_is}^{max} \coloneqq \frac{\Delta_E^s}{\alpha_s}, \delta_{v_is}^{min} \coloneqq \frac{\Delta_{E_{low}}^s}{\alpha_s}$.
8:        If $\delta_{v_is} < \delta_{v_is}^{min}$, GO TO line 4.
9:        Define $\varepsilon_{route}^i \coloneqq [A_i + t_{v_is} + \mu_s, A_i + t_{v_is} + \mu_s + \delta_{v_is}]$.
10:       **for** $\varepsilon_{route}^i \in$ available time intervals of $s$, with time intervals sorted by their starting time
11:       Compute the duration $\delta$ of $\varepsilon_{route}^i \cap \varepsilon_s$ if $\varepsilon_{route}^i \cap \varepsilon_s \neq \emptyset$; otherwise, GO TO line 2.
12:       **if** $\delta \geq \delta_{v_is}^{max}$
13:         Insert recharge at $i$ with charging time $\delta_{v_is}^{max}$ and add this charging event to $s$.
14:         set $\Delta_E \coloneqq 0$ and **return true**.
15:       **else if** $\delta_{v_is}^{min} \leq \delta < \delta_{v_is}^{max}$
16:         Insert recharge at $i$ with charging time $\delta$ and add this charging event to $s$.
17:         Update $\Delta_E \coloneqq \Delta_E^s - \delta \times \alpha_s$ and update the $E_j$ with position $j$ after $i$.

| | |
|---|---|
| 18: | If $\Delta_E > 0$, update $i_{low}$ and $\Delta_{E_{min}}$ to the next position with its SoC below $E_{min}$, and GO TO line 2; otherwise, **return true**. |
| 19: | **else** |
| 20: | GOTO line 10. |
| 21: | **end if** |
| 22: | **end for** |
| 23: | **end for** |
| 24: | **end for** |
| 25: | **return false** and remove the inserted charging events if $\Delta_E > 0$. |

## 5  Computational experiments

To evaluate the proposed algorithm, we generate a set of test instances for the computation experiments. The hybrid LNS algorithm is implemented in the Julia programming language. All experiments are conducted on a laptop equipped with an 11th Gen Intel(R) Core(TM) i5-1135G7 CPU and 64 GB of RAM. The test instance generation is first described, followed by a sensitivity analysis with respect to the key parameters of the algorithm. The performance of destroy, repair, and local search operators designed specifically for the EIDARP is analyzed. The proposed hybrid LNS is compared with the solution found by the commercial solver Gurobi v11.0, given an 8-hour computational time limit.

### 5.1  Test instances

Test instances are generated based on a scenario with two crossed transit lines in a 16 by 16 km square area (see Figure 3). The two transit lines operate bi-directionally, and each has three transit stations with one transfer station in the middle, allowing customers to transfer from one transit line to another. The EIDARP service planning horizon is 2 hours, with a total of four departures (two for each direction) operating on each line. The average operating speed of the transit is 50 km/h, while the bus speed is set as 25 km/h. The bus operator operates two types of buses, with two depots and a sufficiently large fleet size, to serve its customers. Each bus must start and end its service at its designated depot. Each depot is equipped with an operator-owned DC fast charger. We assume a service time (1 minute) for each charging operation. The parameter setting is presented in Table 2. Note that the penalty cost $\omega$ for customer rejection is set to a larger value, ensuring the system serves as many customers as possible. Five test instances are generated for the algorithm parameter tuning, each consisting of 50 customers with origins and destinations randomly located in the service area. Each instance consists of an equal number of outbound and inbound customers. The performance of the algorithm is also evaluated based on a set of test instances with the number of customers ranging from 10 to 100.

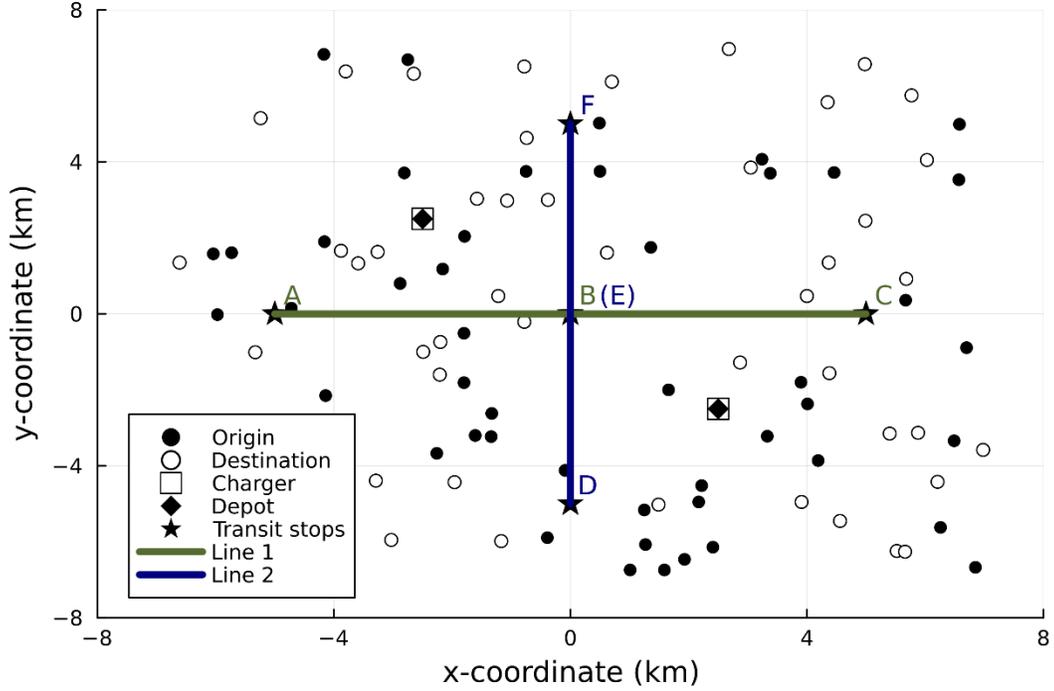

**Figure 3:** Example of the test scenario with two bi-directional transit lines and 50 customers. Line 1 is A-B-C, and line 2 is D-E-F, where B and E are transfer stations.

**Table 2**: The parameter settings of the numerical study.

| Parameter | Value |
| --- | --- |
| Planning horizon $T$ | 2 hours |
| Number of transit lines and stations | 2 lines and 5 stations |
| Number of transit departures | 4 per line (2 per direction) |
| Number of depots | 2 |
| Number of chargers | 2 |
| Number of bus types | 2 types |
| Bus capacity $Q^k$ | 15 and 22 passengers/bus |
| Bus battery capacity | 69 and 103.5kWh |
| Bus energy consumption rate | 0.552 and 0.828 kWh/km |
| Bus average speed | 25km/h for both types |
| Maximum (Minimum) SoCs of buses | 80% (10%) of the battery capacity |
| Charging rate | 0.83kWh/minute (DC fast charger) |
| Time window width at the customer's origin or destination | 15 minutes |
| Maximum waiting time at transit stations $\gamma$ | 10 minutes |
| Minimum and maximum transfer time $\eta^{min}, \eta^{max}$ | 1 minute, 10 minutes |
| Maximum walking distance | 1.5 km |
| Walking speed | 0.085 km/minute |
| Detour factor $\varphi$ | 1.5 |
| The weights in the objective function ($\lambda_1$ and $\lambda_2$) | 1, 1 |
| Penalty cost $\omega$ | 200 |
| Service time $\mu$ | 0.5 minute |
| Service time for each charging operation | 1 minute |

## 5.2 Algorithm parameter tuning

In this section, we configure five parameters, i.e., $n_{iter}, t_{max}, T_{red}, \xi_{max}$, and $\alpha$ for the proposed hybrid LNS. Based on the preliminary experiments, these parameters are initially set to $(n_{iter}, t_{max}, T_{red}, \xi_{max}, \alpha) = (500, 1.1, 100, 0.2, 1.02)$. Each parameter is then varied within a range around its initial value to assess its impact on performance. The experiments are conducted sequentially as the parameters are independent of each other, except for $t_{max}$ and $T_{red}$. Each time a new value for the analyzed parameter is identified, it is updated in subsequent analyses. Note that more advanced parameter tuning approaches can also be used (Lindauer et al., 2022; López-Ibáñez et al., 2016). Five instances with 50 customers are applied for different parameter values, and the computational (CPU) time and objective function values are recorded after five runs. The best-known solution (BKS) of each instance is recorded from all the conducted experiments. The average gap reported in the results is calculated as the average gap compared with the 8-hour BKS obtained by Gurobi. The results are presented in Tables 3–5.

First, $n_{iter}$ is analyzed with the results presented in Table 3. The solution quality increases with the number of iterations at the cost of higher CPU time. After 600 iterations, the gap reduces slightly, but CPU time increases by more than 20%. Therefore, $n_{iter}$ is set to 600 by considering the trade-offs of additional computational time and the reduction in the objective function.

**Table 3:** Results of sensitivity analysis on $n_{iter}$.

| $n_{iter}$ | 300 | 400 | 500 | 600 | 700 | 800 | 900 | 1000 |
|---|---|---|---|---|---|---|---|---|
| Avg. Gap (%) | 2.36 | 2.13 | 2.03 | 1.73 | 1.68 | 1.67 | 1.65 | 1.61 |
| Avg. CPU time (s) | 166 | 258 | 327 | 397 | 481 | 492 | 548 | 623 |

As $t_{max}$ and $T_{red}$ are used to determine the DA acceptance criteria jointly, the experiments are conducted by setting 100 combinations with 10 values each for both $t_{max}$ and $T_{red}$. Table 4 presents the average gap of each value on these two parameters. Among all possible combinations, the average gap is smallest when $t_{max} = 1.1$ and $T_{red} = 700$. As the CPU time does not exhibit a clear trend across the tested combinations, we set $t_{max} = 1.1$ and $T_{red} = 700$.

**Table 4:** Results of sensitivity analysis on $t_{max}$ and $T_{red}$.

| $t_{max}$ | 0.7 | 0.8 | 0.9 | 1.0 | 1.1 | 1.2 | 1.3 | 1.4 | 1.5 | 1.6 |
|---|---|---|---|---|---|---|---|---|---|---|
| Avg. gap (%) | 2.05 | 1.86 | 2.16 | 1.63 | 1.76 | 1.81 | 1.93 | 1.82 | 1.89 | 2.06 |
| $T_{red}$ | 100 | 200 | 300 | 400 | 500 | 600 | 700 | 800 | 900 | 1000 |
| Avg. gap (%) | 1.96 | 2.10 | 2.00 | 1.90 | 1.83 | 1.81 | 1.73 | 1.85 | 2.02 | 1.77 |

Note: The minimum gap is 1.60% for $t_{max} = 1.1$ and $T_{red} = 700$.

Table 5 presents the results for $\xi_{max}$, $\alpha$, and $\beta$. $\xi_{max}$ is the parameter of the degree of destruction. We can observe that when $\xi_{max}$ increases, the CPU time increases, but the average gap does not increase accordingly. The minimum average gap is 1.73% with $\xi_{max} = 0.25$. $\alpha$ controls the frequency of executing the local search procedure. The higher the $\alpha$ is, the longer the CPU time and the better the average gap. We set $\alpha$ to 1.06 at the minimum average gap.

**Table 5:** Results of sensitivity analysis on $\xi_{max}$ and $\alpha$.

| $\xi_{max}$ | 0.15 | 0.2 | 0.25 | 0.3 | 0.35 |
|---|---|---|---|---|---|
| Avg. gap (%) | 1.94 | 1.77 | 1.73 | 1.77 | 1.92 |
| Avg. CPU time (s) | 253 | 343 | 459 | 555 | 681 |
| $\alpha$ | 1.04 | 1.05 | 1.06 | 1.07 | 1.08 |
| Avg. gap (%) | 1.94 | 1.64 | 1.43 | 1.45 | 1.80 |
| Avg. CPU time (s) | 488 | 509 | 525 | 561 | 590 |

In summary, the best configuration for the examined parameters is $(n_{iter}, t_{max}, T_{red}, \xi_{max}, \alpha) = (600, 1.0, 700, 0.25, 1.06)$, which is utilized in the following experiments.

### 5.3 Performance of the algorithm operators

This section evaluates the six operators tailored for EIDARP: $r_{TS}, r_{Tpriority}, d_{related}$, and three local search operators: *exchange TPs of customers, replace TP by bus*, and *re-assign first-/last-mile bus service*. To evaluate their impact on the solution quality, each is removed from the hybrid LNS and compared to the case in which all operators are included. The experiment is conducted on the five 50-customer test instances, each solved with five runs.

Table 6 presents the average gap to BKS and CPU time for the different cases. The average gap, excluding one of these operators, results in 0.41% to 0.84% higher than when all operators are included. In the case without $r_{Tpriority}$, the average gap increases up to 2.09% with a slight reduction in the average CPU time, indicating $r_{Tpriority}$ is an effective and efficient operator. By contrast, removing the local search operator *exchange TPs of customers* results in the smallest improvement to the average gap, but not the greatest reduction in CPU time. This suggests that this operator has a relatively lower impact compared with the other five. Despite these variations, the six operators demonstrate their influence on the hybrid LNS algorithm for EIDARP. The exclusion of any operator negatively impacts the solution quality.

**Table 6:** Performance of the EIDARP-excusive operators.

| | Avg. Gap | Avg. CPU time (s) |
|---|---|---|
| **All operators** | **1.43%** | **525** |
| Without $r_{TS}$ | 2.14% | 346 |
| Without $r_{Tpriority}$ | 2.09% | 501 |
| Without $d_{related}$ | 1.95% | 418 |
| Without local search operator *exchange TPs of customers* | 1.84% | 459 |
| Without local search operator *replace TP by bus* | 1.85% | 522 |
| Without local search operator *re-assign first-/last-mile bus service* | 2.27% | 444 |

### 5.4 Algorithm performance

The solution quality of the proposed hybrid LNS is compared with the solutions of the MILP model solved by Gurobi. Test instances are generated using the scenario described in Section 5.1, with the number of customers ranging from 10 to 50 when buses are fully charged at the beginning of the service (full initial SoCs of buses), and with 10 to 20 customers when the initial SoCs of buses are set to 30%. The results from Gurobi are obtained after 8 hours of CPU time, and the results of the hybrid LNS are based on five runs. The test instances are available in github.com/YMF2022/EIDARP-instances.

Table 7 summarizes the results when buses are fully charged at the beginning of the service. Among the instances, Gurobi obtained the optimal solution only for the 10-customer case with 1239 seconds; the average gap to the lower bound is 32.1% (see the MILP column). The hybrid LNS obtains the exact solution for the 10-customer test instance with an average of 6 seconds of CPU time. For the remaining instances, the best objective function value obtained by the hybrid LNS over the five runs consistently outperforms the MILP solutions obtained by Gurobi, with an average improvement of -23.8% (best) and -23.3% (average). The 8-hour Gurobi solution has an average gap of 32.1% to its lower bound. In contrast, our hybrid LNS algorithm achieves an average gap of 12.9% to the Gurobi lower bound and improves the objective value of the Gurobi solution by 23.3 %. This suggests a good quality of the hybrid LNS solution. Moreover, the hybrid LNS obtains the solutions significantly faster, with an average of 136 seconds only.

**Table 7:** Performance of hybrid LNS and its comparison with MILP (full initial SoC of buses).

| $n$ | Gurobi (MILP) | | Hybrid LNS | | | | |
|---|---|---|---|---|---|---|---|
| | Obj. value[1] | Gap to LB[2] | Best obj. | Best gap[3] | Avg. gap[3] | Best gap to LB[4] | CPU avg. (s) |
| 10 | 490.82 | 0.0%[5] | 490.82 | 0.0% | 0.0% | 0.0% | 6 |
| 15 | 688.48 | 9.5% | 683.52 | -0.7% | 0.0% | 8.9% | 12 |
| 20 | 815.54 | 2.4% | 811.69 | -0.5% | 0.5% | 1.9% | 11 |
| 25 | 1203.21 | 22.3% | 1072.10 | -10.9% | -10.8% | 12.8% | 33 |
| 30 | 1419.42 | 24.7% | 1197.97 | -15.6% | -14.3% | 10.8% | 76 |
| 35 | 2640.80 | 51.8% | 1522.13 | -42.4% | -42.1% | 16.4% | 137 |
| 40 | 3351.43 | 63.5% | 1759.03 | -47.5% | -47.0% | 30.5% | 162 |
| 45 | 3457.61 | 57.9% | 1792.33 | -48.2% | -47.8% | 18.7% | 365 |
| 50 | 3727.94 | 56.9% | 1917.23 | -48.6% | -48.3% | 16.2% | 417 |
| Average | 1977.25 | 32.1% | 1249.65 | -23.8% | -23.3% | 12.9% | 136 |

Notes: 1. Objective function value of the incumbent solution after 8-hour CPU time. 2. Gap to the lower bound. 3. Gap compared with the objective function value found by Gurobi. 4. Gap between the best objective value of the hybrid LNS and the lower bound of the Gurobi solution. 5. Optimal solution with CPU time 1239 seconds. Note that no customers are rejected.

In the low initial SoC (30% of the battery capacity) scenario, the gaps to the lower bound of MILP are significantly higher than in the fully initial SoC scenario. This experiment limits the maximum number of customers to 20, as Gurobi was unable to obtain feasible solutions, or the obtained gaps to the lower bound are relatively high for larger problem sizes (a gap of 78.0% is observed for a 30-customer scenario). The results are provided in Table 8. The instance with 10 customers is the only one for which MILP finds the optimal solution. However, due to recharging scheduling in the low SoC case, the CPU time is 7.3 hours, which is 20 times longer than the fully charged scenario presented in Table 7. By contrast, the hybrid LNS shows a 1-second difference for the same instance. The MILP exhibits an average gap of 12.5 % to its lower bound across all instances. In contrast, the best solutions found by the hybrid LNS reduce this gap to just 5.5 % relative to the same bound. When comparing directly to the MILP's objective values, the hybrid LNS achieves a best-case improvement of 5.5 % and an average improvement of 4.8 %. When buses have low initial SoCs, the variance in the hybrid LNS performance is slightly higher than in the fully charged case. In terms of charging times, the hybrid LNS achieved an average charging time of 7.9 minutes, compared with MILP's average of 17.4 minutes. This difference aligns with the hybrid LNS strategy of recharging only the minimum required amount, as described in Section 4.6.

**Table 8:** Performance of the hybrid LNS and its comparison with MILP (low initial SoC of buses).

| $n$ | Gurobi (MILP) | | | Hybrid LNS | | | | | |
|---|---|---|---|---|---|---|---|---|---|
| | Obj. value[1] | Gap to LB[2] | Charging time (min) | Best obj. | Best gap[3] | Avg. gap[3] | Best gap to LB[4] | Charging time (min) | Avg. CPU (s) |
| 10 | 490.82 | 0.0%[5] | 5.0 | 490.82 | 0.0% | 1.6% | 0.0% | 5.0 | 7 |
| 11 | 502.37 | 4.1% | 7.8 | 502.37 | 0.0% | 0.2% | 3.9% | 6.1 | 11 |
| 12 | 584.67 | 15.1% | 11.0 | 573.31 | -1.9% | -0.5% | 11.5% | 13.7 | 8 |
| 13 | 486.89 | 2.8% | 16.3 | 486.89 | 0.0% | 1.6% | 2.1% | 1.7 | 6 |
| 14 | 653.80 | 11.7% | 0.1 | 610.56 | -6.6% | -6.5% | 4.1% | 0.0 | 5 |
| 15 | 790.14 | 27.2% | 4.7 | 689.35 | -13.5% | -12.9% | 9.1% | 2.9 | 16 |
| 16 | 685.22 | 14.8% | 8.2 | 650.39 | -6.3% | -5.6% | 7.0% | 7.8 | 12 |
| 17 | 889.40 | 20.2% | 36.1 | 808.23 | -9.7% | -9.3% | 7.9% | 10.4 | 25 |
| 18 | 686.44 | 8.7% | 4.1 | 675.83 | -1.5% | -1.5% | 6.5% | 4.3 | 15 |
| 19 | 998.44 | 23.6% | 1.8 | 851.03 | -14.8% | -14.7% | 5.1% | 6.0 | 19 |
| 20 | 862.16 | 10.4% | 96.0 | 827.77 | -5.8% | -4.8% | 3.8% | 29.3 | 13 |
| Average | 693.67 | 12.5% | 17.4 | 653.04 | -5.5% | -4.8% | 5.5% | 7.9 | 13 |

Notes: 1. Objective function value of the incumbent solution after 8-hour CPU time. 2. Gap to the lower bound. 3. Gap compared with the objective function value found by Gurobi. 4. Gap between the best objective value of the hybrid LNS and the lower bound of Gurobi's 8-hour solution. 5. Optimal solution with CPU time 26118 seconds. Note that no customers are rejected.

Table 9 presents the results of the hybrid LNS for the test instances with up to 100 customers, featuring both full and low initial SoCs of buses. In both scenarios, the CPU time increases significantly as the number of customers increases. No recharging events are observed in the full SoC scenario, whereas the low SoC scenario shows an overall increase in charging times as the number of customers grows. Although the low SoC scenario consistently exhibits longer charging times than the full SoC scenario, the difference in their CPU times is minor, suggesting that the recharging scheduling approach is efficient and does not significantly impact the CPU time for these experiments. Besides, in the low SoC instance with 100 customers, the nine-step feasibility evaluation scheme accounts for 2% of total CPU time, while the recharging scheduling algorithm accounts for 1 %. The CPU time is mainly occupied by the repair procedures (71%).

**Table 9:** Performance of the hybrid LNS on full and low initial SoC.

| $n$ | Full (100%) initial SoC | | | Low (30%) initial SoC | | | |
|---|---|---|---|---|---|---|---|
| | Best obj. | Avg. obj. | Avg. CPU (s) | Best obj. | Avg. obj. | Avg. charging time (min) | Avg. CPU (s) |
| 10 | 490.82 | 490.82 | 10 | 490.82 | 498.90 | 5.0 | 7 |
| 20 | 811.69 | 819.58 | 11 | 811.84 | 820.62 | 7.9 | 13 |
| 30 | 1197.97 | 1216.87 | 76 | 1211.10 | 1219.12 | 10.6 | 134 |
| 40 | 1759.03 | 1777.56 | 162 | 1787.43 | 1801.12 | 13.2 | 251 |
| 50 | 1917.23 | 1928.69 | 417 | 1953.05 | 1976.15 | 22.7 | 513 |
| 60 | 2467.26 | 2492.48 | 1031 | 2499.64 | 2525.25 | 26.0 | 1125 |
| 70 | 2630.24 | 2638.84 | 1078 | 2697.17 | 2715.59 | 14.8 | 1361 |
| 80 | 3060.91 | 3122.02 | 1821 | 3137.84 | 3158.55 | 20.7 | 2891 |
| 90 | 3256.01 | 3267.10 | 3722 | 3312.57 | 3350.95 | 21.9 | 4299 |
| 100 | 3864.78 | 3908.28 | 6163 | 4007.95 | 4029.80 | 27.4 | 7054 |
| Average | 2145.59 | 2166.22 | 1449 | 2575.00 | 2595.03 | 17.0 | 1765 |

Figure 4 (a) illustrates the average CPU time for the full and low SoC scenarios, confirming that CPU time increases exponentially with the increase in the number of customers. The charging scheduling part for the low SoC scenario does not have a significant impact on its CPU time. Figure 4 (b) also compares the CPU time between EDARP and EIDARP under a full SoC scenario, where no recharging is observed for either case. A significantly higher CPU time is observed for solving the EIDARP compared with that of EDARP. The increase in CPU time for EIDARP is attributed to the additional complexity of considering different possible travel options, including walk, bus, and transit (travel options 1-4). EDARP handles only one travel option (direct bus trips from customer origins to destinations).

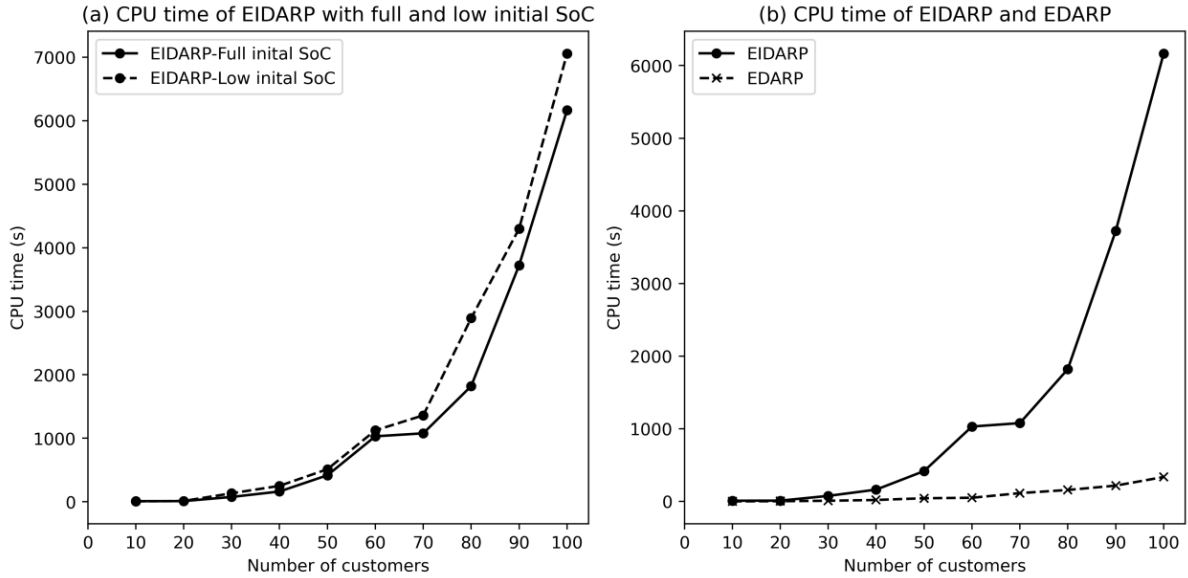

**Figure 4:** Computational time of EIDARP at different initial SoC levels and comparison with EDARP.

## 6   Analyzing the performance of the EIDARP

The integrated service leverages existing fixed transit services, and its quality is influenced not only by bus operators but also by the operations of the transit system as a whole. Molenbruch et al. (2021) investigated the service performance by examining different configurations for demand-responsive bus operators (e.g., maximum detour) and transit operators (e.g., transit speed and frequency). Their findings indicated that the transit frequency is less critical due to the flexibility provided by demand-responsive vehicles.

This section investigates service performance based on the characteristics of the proposed EIDARP from both the bus operator's and customers' perspectives. In the problem, we identify seven key parameters for the performance analysis, which are:
(1) detour factor $\varphi$,
(2) weights on customer travel time in the objective function $\lambda_2$,

(3) maximum waiting time $\gamma$,
(4) bus fleet size,
(5) bus operation speed,
(6) initial SoC of buses,
(7) transit networks.

The analysis is based on a baseline scenario of 100 customers with the problem configuration described in Section 5.1. The reference values of the EIDRT system parameters are set as $\varphi = 1.5$, $\lambda_2 = 1.0$, unlimited bus fleet size, 25km/h for the average bus speed, full (100%) initial SoC of buses, and two crossed transit lines. The subsequent experiments modify only one parameter at a time to assess the service performance. The following key performance indicators (KPIs) are identified for evaluation.
- BTT: <u>Total</u> bus travel time in minutes
- RT: <u>Total</u> recharging time in minutes
- CTT: <u>Average</u> customer travel time: in-transit time, in-bus time plus walking time in minutes
- CTT-transit: <u>Average</u> customer in-transit time in minutes
- CTT-bus: <u>Average</u> customer in-bus time in minutes
- CTT-walk: <u>Average</u> customer walking time in minutes
- WT: <u>Average</u> customer experienced waiting time in minutes, including waiting time at transit station for both first and last mile, along with waiting time during bus services.
- # cus-transit: Number of customers using the transit service.
- # used buses: Number of used buses
- # cus/bus: <u>Average</u> number of served customers per bus.
- # reject: Number of rejected customers
- # TO4 (bus-transit-bus): Number of customers using travel option 4 (bus + transit + bus). The number of customers using the other options can be calculated as follows.
  - The number of customers using travel options 1 to 3 can be obtained by the number of customers using transit service minus # TO4.
  - The number of customers using travel option 5 = the total number of customers minus the total number of customers using the transit service.

## 6.1 Detour factor $\varphi$

This section analyzes the effects of the detour factor, $\varphi$, within the range of 1.3 to 2.5. The results are summarized in Table 10. The result shows that an increase in $\varphi$ leads to higher CTT and reduced BTT. The increase in CTT can be attributed to higher utilization of the transit service and longer waiting times, as indicated by the rise in CTT-transit and WT. Among customers using transit services, most use travel option 4 (see column #TO4), and a few walk to transit stations (few WT). However, the values of CTT-bus increase slightly with larger $\varphi$ due to a higher degree of ride-sharing, evidenced by a decrease in the number of used buses and an increase in the average number of customers served per bus. WT also increases accordingly, reaching up to 6.4 minutes when $\varphi = 2.5$, while no consistent trend is observed for CTT-walk.

**Table 10:** Impacts of $\varphi$ on the performance of EIDARP.

| $\varphi$ | BTT | # used buses | # cus/bus | # cus-transit | CTT | CTT-transit | CTT-bus | CTT-walk | WT | #TO4 |
|---|---|---|---|---|---|---|---|---|---|---|
| 1.3 | 2232.0 | 36 | 3.5 | 31 | 20.5 | 1.9 | 18.6 | 0.00 | 1.9 | 31 |
| 1.5[1] | 1902.9 | 29 | 4.6 | 39 | 21.4 | 2.8 | 18.4 | 0.15 | 2.7 | 36 |

| | | | | | | | | | |
|---|---|---|---|---|---|---|---|---|---|
| 1.7 | 1864.0 | 28 | 5.1 | 52 | 22.3 | 4.0 | 18.1 | 0.24 | 3.9 | 49 |
| 1.9 | 1829.0 | 29 | 5.0 | 53 | 23.0 | 3.7 | 19.2 | 0.14 | 4.5 | 51 |
| 2.1 | 1770.5 | 26 | 5.7 | 53 | 23.6 | 3.9 | 19.6 | 0.17 | 4.8 | 51 |
| 2.3 | 1783.4 | 26 | 5.4 | 47 | 24.2 | 3.9 | 20.3 | 0.07 | 5.0 | 46 |
| 2.5 | 1775.1 | 25 | 6.0 | 59 | 24.3 | 4.6 | 19.5 | 0.25 | 6.4 | 56 |

Note: 1. The baseline scenario. No customers were rejected, and no recharging operations were required for all the cases.

Figure 5 illustrates the detour experienced by customers. Figure 5 (a) shows a more dispersed distribution of customers' experienced detours when a larger $\varphi$ is used. For each value of $\varphi$, some customers experience the upper bound of maximum travel time, whereas the minimum detours experienced by customers remain at a similar level, regardless of the different values of $\varphi$. When customers' experienced detours are smaller than 1.0, they benefit from using transit service (operating with a higher speed (50 km/h) compared with the bus (25 km/h)), as shown in Figure 5 (b). Figure 5(b) classifies customers into two groups: those (partially or exclusively) served by transit (travel options 1-4) and those served exclusively by buses (travel option 5). As presented in Figure 5 (b), customers served by transit could eventually experience a detour less than 1.0, while customers exclusively served by bus always experience a detour equal to or greater than 1.0. However, customers with travel options 1-4 reach the upper bound of the maximum detour more often than those served by bus only.

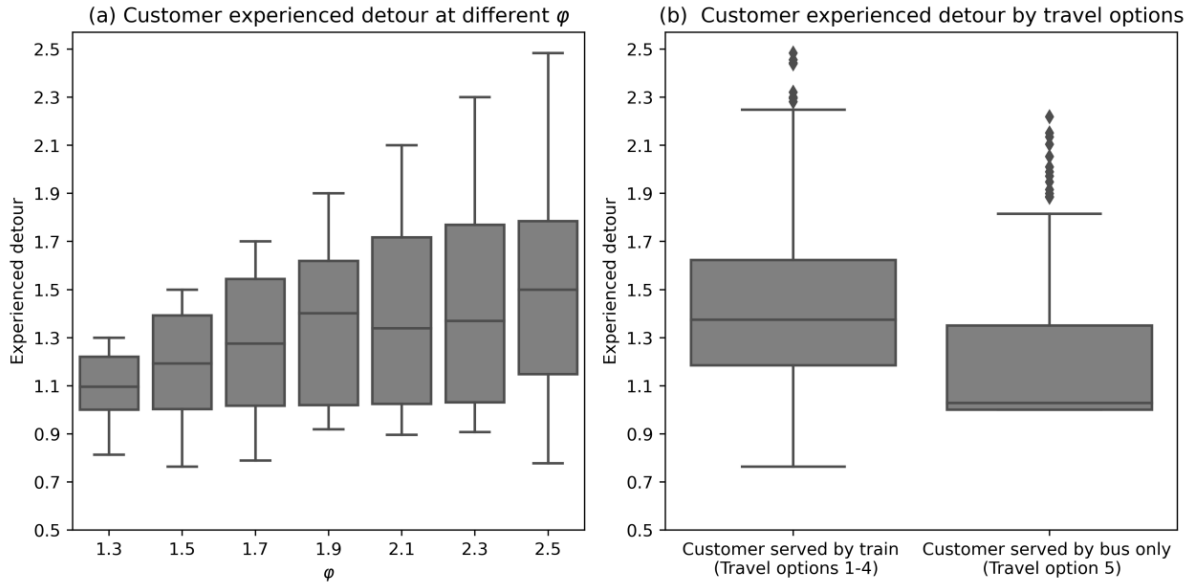

**Figure 5:** Impacts of $\varphi$ on customers' experienced detour.

## 6.2 Weight on customer travel time in the objective function $\lambda_2$

Recall that $\lambda_1, \lambda_2,$ and $\omega$ are the weights for bus travel time, customer travel time, and penalty cost of rejecting customers used in the objective function, respectively. The experiments vary $\lambda_2$ from 0.0 to 2.5 while keeping $\lambda_1 = 1.0$ and $\omega = 1.0$. The results are presented in Table 11. As $\lambda_2$ increases, CTT decreases while BTT increases. The increase in BTT is due to more customers utilizing bus service, but fewer using transit service, as reflected in the changes in CTT-transit and CTT-bus. Moreover, when more customers are served by door-to-door DRT buses, a decrease in CTT-walk is observed. Unlike other KPIs, WT does not exhibit a consistent trend.

**Table 11:** Impacts of $\lambda_2$ on the performance of EIDARP.

| $\lambda_2$ | BTT | # used buses | # cus/bus | # cus-transit | CTT | CTT-transit | CTT-bus | CTT-walk | WT | # TO4 |
|---|---|---|---|---|---|---|---|---|---|---|
| 0.0 | 1884.8 | 31 | 4.2 | 34 | 23.1 | 2.8 | 20.0 | 0.25 | 1.5 | 30 |
| 0.5 | 1893.9 | 30 | 4.5 | 38 | 21.5 | 2.9 | 18.6 | 0.00 | 3.0 | 38 |
| 1.0[1] | 1902.9 | 29 | 4.6 | 39 | 21.4 | 2.8 | 18.4 | 0.15 | 2.7 | 36 |
| 1.5 | 1934.0 | 29 | 4.4 | 32 | 20.5 | 2.3 | 18.1 | 0.08 | 2.8 | 30 |
| 2.0 | 1986.1 | 31 | 4.4 | 37 | 19.9 | 2.3 | 17.6 | 0.00 | 2.9 | 37 |
| 2.5 | 2055.8 | 33 | 3.9 | 35 | 19.7 | 2.2 | 17.5 | 0.00 | 1.8 | 35 |

Note: 1. The baseline scenario. No customers rejected and no recharging operations for all the cases.

As the customer travel time in the objective function of the EIDARP excludes waiting time, while the maximum travel time constraint accounts for it, both experienced detour and experienced detour without waiting time are defined.

- Experienced detour for $r \in R$: Customer $r$'s travel time plus waiting time divided by their respective direct travel time by bus.
- Experienced detour without WT for $r \in R$: Customer $r$'s travel time divided by their respective direct travel time by bus.

Figure 6 investigates the impact of $\lambda_2$ on customers' experienced detour with (Figure 6(b)) and without waiting times (Figure 6(a)). As the detour factor $\varphi$ is set to 1.5, both KPIs are less than or equal to this value. In Figure 6 (a), an increase in $\lambda_2$ corresponds to a reduction in the experienced detour without WT, with the distribution becoming more concentrated around 1.0. This trend reflects a higher proportion of customers served by direct bus services without shared rides, which is also observed in the decrease of #cus/bus in Table 11. Similarly, customers' experienced detour decreases with the increase in $\lambda_2$, though not as significantly as observed in Figure 6(a). This is because the objective function minimizes customer travel times (without WT), leading to a faster reduction in customers' experienced detours without WT. Besides, customers' experienced detour still reaches the upper bound of 1.5. This is attributed to the waiting time not decreasing proportionally, as indicated in Table 11.

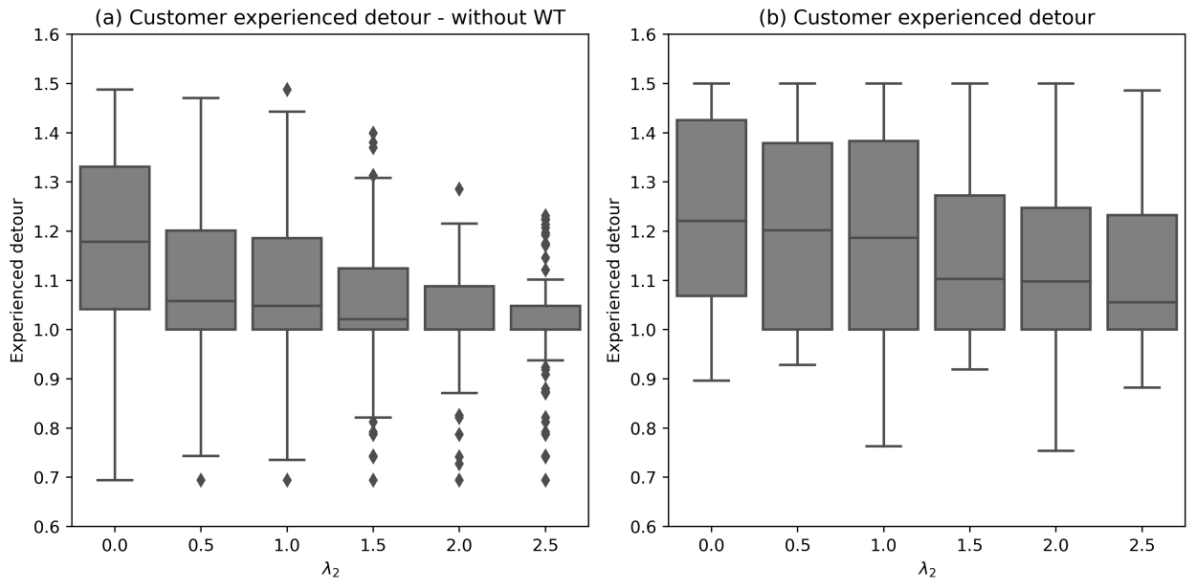

**Figure 6:** Impacts of $\lambda_2$ on customers' experienced detour without WT and experienced detour.

## 6.3 Maximum waiting time $\gamma$

This section investigates the maximum waiting time constraint at transit stations when customers transfer between bus and transit services. $\gamma$ is varied from 5 to 15 minutes to evaluate customers' experienced waiting times. Customers' experienced waiting time includes waiting time at transit stations for both first- and last-mile connections, as well as the waiting time during the bus service. In Figure 7(a), the customer waiting time at each transit station is constrained when $\gamma = 5$ or $\gamma = 6$, but when $\gamma > 7$, the customer's waiting time ranges between 0 and 14.5 minutes, regardless of the value of $\gamma$. This is because customers' experienced waiting time is partially bounded by the maximum travel time constraint. The variability in Figure 7(a) is primarily attributed to the customers served by transit, as depicted in Figure 7(b). The group of customers served by transit experiences more diverse and longer waiting times, whereas the bus-only group experiences shorter and more consistent waiting times, typically less than 2 minutes.

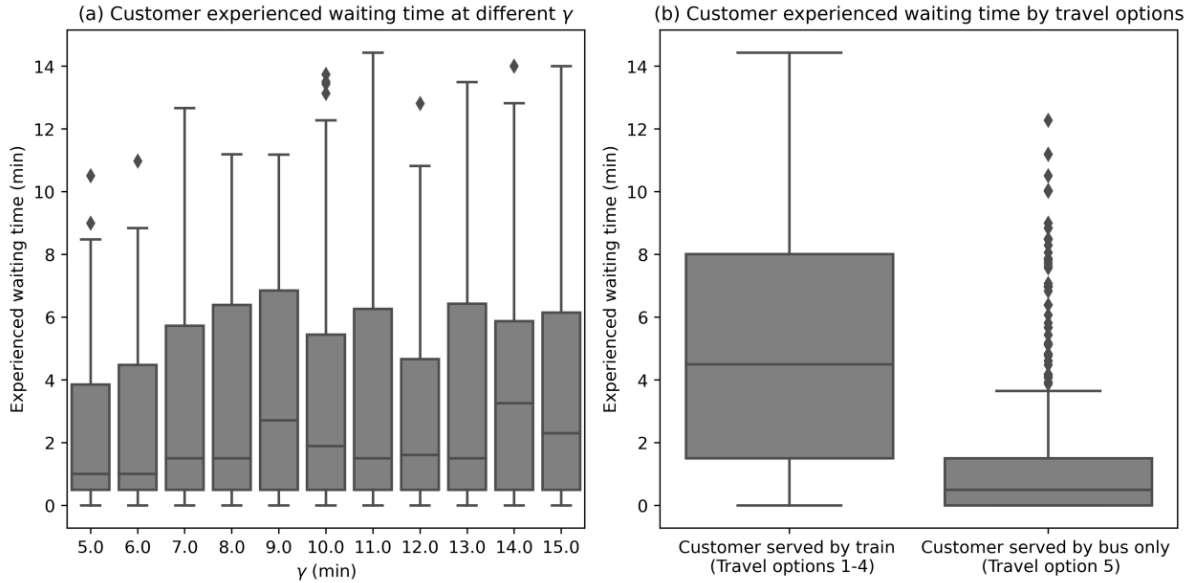

**Figure 7:** Impacts of $\gamma$ on customers' experienced waiting time.

## 6.4 Bus fleet

Table 12 shows the impacts of fleet sizes. When the fleet size is 10 (the minimum value), the number of rejected customers reaches 52, and no customers use the transit service. As the fleet size increases, customer rejections decrease, and more BTT is required to serve additional customers. In the case of an unlimited fleet size, all customers can be served with 29 buses, but the CTT increases. Both CTT-transit and CTT-bus show little change with increasing fleet size. However, WT increases from 0.4 minutes (fleet size is 10 buses) to 2.7 minutes (fleet size is 29 buses). This increase is likely due to longer waiting times at transit stations and the added waiting time for customers sharing rides. The average number of served customers per bus does not exhibit a consistent trend.

**Table 12:** Impacts of bus fleet size on the performance of EIDARP.

| Fleet size | BTT | # used buses | # cus/bus | # reject | # cus-transit | CTT | CTT-transit | CTT-bus | CTT-walk | WT | # TO4 |
|---|---|---|---|---|---|---|---|---|---|---|---|
| 10 | 780.8 | 10 | 4.8 | 52 | 0 | 15.7 | 0.0 | 15.7 | 0.00 | 0.4 | 0 |
| 15 | 1223.2 | 15 | 5.6 | 31 | 17 | 18.1 | 1.8 | 16.2 | 0.10 | 1.8 | 16 |

| 20 | 1504.7 | 20 | 5.5 | 16 | 30 | 20.4 | 2.5 | 17.8 | 0.08 | 2.3 | 29 |
| 25 | 1868.2 | 25 | 4.8 | 3 | 26 | 21.6 | 1.9 | 19.7 | 0.00 | 2.1 | 26 |
| ULD[1] | 1902.9 | 29 | 4.6 | 0 | 39 | 21.4 | 2.0 | 18.4 | 0.15 | 2.7 | 36 |

Note: 1. ULD means the unlimited fleet size, which is the baseline scenario. No recharging operations for all cases.

## 6.5 Bus operational speed

As different service areas (rural and urban) limit bus travel speeds, this section examines the impact of varying the average bus speed from 25 to 50 km/h. The results, presented in Table 13, show that increasing bus speed reduces both BTT and CTT. The number of served customers per bus increases, resulting in a reduction in the total number of buses required. When bus speed doubles from 25 km/h to 50 km/h, the BTT is reduced by more than half due to an increase in shared rides. When the bus speed reaches 50 km/h (the same as transit operational speed), only six customers utilize the transit service. This observation highlights that increasing the bus speed undermines the purpose of transit integration.

**Table 13:** The impacts of bus travel speed on the performance of the EIDARP.

| Speed (km/h) | BTT | # used buses | # cus/bus | # cus-transit | CTT | CTT-transit | CTT-bus | CTT-walk | WT | # TO4 |
|---|---|---|---|---|---|---|---|---|---|---|
| 25[1] | 2047.7 | 32 | 4.3 | 42 | 21.3 | 3.1 | 18.1 | 0.11 | 3.0 | 40 |
| 30 | 1535.4 | 26 | 4.6 | 23 | 17.6 | 1.7 | 15.9 | 0.00 | 2.2 | 23 |
| 35 | 1256.3 | 22 | 5.4 | 18 | 15.8 | 1.1 | 14.7 | 0.00 | 1.3 | 18 |
| 40 | 1118.1 | 20 | 5.9 | 18 | 13.7 | 1.0 | 12.7 | 0.00 | 1.1 | 18 |
| 45 | 986.5 | 19 | 5.7 | 10 | 12.2 | 0.6 | 11.6 | 0.01 | 1.0 | 9 |
| 50 | 868.2 | 18 | 5.9 | 6 | 11.0 | 0.4 | 10.7 | 0.00 | 0.7 | 6 |

Note: 1. The baseline scenario. No customers were rejected, and no recharging operations were required for all the cases.

## 6.6 Initial SoC of buses

The previous results do not examine bus recharging operations when buses are fully recharged at the beginning of the service. However, in reality, there may be some uncertainties, meaning that a 100% initial SoC cannot be ensured for all buses. This section analyzes the impact of recharging on this service by setting an identical initial SoC level for all buses, ranging from 20% to 100%. The results are presented in Table 14. Even with an unlimited bus fleet, four customers are still rejected when the initial SoC is 20%. This occurs because charging station congestion prevents some buses from recharging in time to meet customer time window constraints despite the availability of the fleet. In this extreme case, we demonstrate that serving all 96 customers with an initial SoC of 20% would require 63 buses. Besides, the case of 20% initial SoC shows the lowest CTT-bus but the highest values in CTT-transit and number of customers using transit, suggesting that buses serve more first-mile and last-mile services instead of providing door-to-door bus connections between customers' origins and destinations. When the initial buses' SoCs increase, BTT decreases, resulting in both a higher number of customers served per bus and a shorter distance to recharge for CSs. We observe that CTT-bus also increases slightly while CTT-transit reduces slightly, suggesting that buses become more efficient in serving customers as a door-to-door bus service when recharging is not necessary. Increased CTT-transit is accompanied by a rise in CTT-walk, while no significant changes occur in WT across different initial SoC settings. With an initial SoC of 20 %, buses accumulate a total charging time of 130.3 minutes over 33 recharge events. At this low SoC, most buses opt to recharge immediately at their depots at the beginning of the operational time. However, due to a very limited number of chargers, the algorithm requires more effort to find feasible charging schedules.

Our computational results show that for the 20% initial SoC case (with 33 recharge events), the CPU time is approximately 50 minutes longer than that for the 25% SoC scenario (with 22 recharge events).

Table 14: The impacts of initial SoC level on the performance of the EIDARP.

| Initial SoC | BTT | # used buses | # cus/bus | RT | # RE[3] | # cus-transit | CTT | CTT-transit | CTT-bus | CTT-walk | WT | # TO[4] | CPU (s) |
|---|---|---|---|---|---|---|---|---|---|---|---|---|---|
| 20%[2] | 2522.7 | 63 | 2.2 | 130.3 | 33 | 46 | 19.9 | 3.9 | 16.0 | 0.00 | 2.7 | 46 | 7057 |
| 25% | 2211.3 | 44 | 3.1 | 61.4 | 22 | 40 | 20.3 | 3.3 | 17.0 | 0.08 | 2.6 | 38 | 4139 |
| 30% | 1989.6 | 33 | 4.0 | 35.9 | 16 | 36 | 20.6 | 2.5 | 18.1 | 0.08 | 2.9 | 34 | 4119 |
| 35% | 1964.3 | 31 | 4.2 | 11.8 | 6 | 34 | 20.8 | 2.7 | 18.2 | 0.01 | 2.5 | 33 | 4021 |
| 40% | 1971.4 | 30 | 4.5 | 1.3 | 1 | 38 | 21.2 | 2.8 | 18.3 | 0.01 | 2.8 | 37 | 4250 |
| 100%[1] | 1902.9 | 29 | 4.6 | 0.0 | 0 | 39 | 21.4 | 2.8 | 18.4 | 0.15 | 2.7 | 36 | 2781 |

Note: 1. The baseline scenario. 2. Four customers are rejected when the initial SoC of buses is set as 20%. No customers are rejected in the other cases. 3. The total number of bus recharge events.

## 6.7 Impact of transit network configurations

This section examines the impact of different transit network configurations on the performance of the hybrid LNS.

### 6.7.1 Transit frequency

We vary the transit headway (i.e., frequency) of both lines from 10 minutes to 60 minutes. The results are shown in Table 15. As transit service is less frequent, both BTT and CTT increase. This is primarily because more customers receive direct, door-to-door bus service, resulting in fewer customers using transit, lower CTT for transit users, and a decrease in the number of customers served per bus. With longer headways, fewer customers use transit, leading to reduced waiting times and shorter walking distances to transit stations. However, the number of used buses shows no clear trend. While increasing transit frequency does reduce both BTT and CTT, the gains are smaller than expected, and the computational time increases exponentially with higher frequencies. Our results are consistent with Molenbruch et al. (2021), who found that the service quality improvements plateau at higher frequencies because the EIDARP's scheduling flexibility already optimizes bus operations to align with scheduled transit services.

Table 15: The impacts of transit frequency.

| Frequency (minute) | BTT | # used buses | # cus/bus | # cus-transit | CTT | CTT-transit | CTT-bus | CTT-walk | WT | # TO[4] | CPU (s) |
|---|---|---|---|---|---|---|---|---|---|---|---|
| 10 | 1840.3 | 30 | 5.0 | 53 | 20.5 | 4.2 | 16.2 | 0.08 | 3.8 | 51 | 34.8[2] |
| 20 | 1764.6 | 30 | 4.7 | 44 | 20.4 | 3.6 | 16.8 | 0.07 | 3.2 | 43 | 3032 |
| 30[1] | 1902.9 | 29 | 4.6 | 39 | 21.4 | 2.8 | 18.4 | 0.15 | 2.7 | 36 | 2781 |
| 40 | 1902.9 | 26 | 4.5 | 18 | 21.7 | 1.2 | 20.5 | 0.00 | 1.6 | 18 | 873 |
| 50 | 1904.6 | 26 | 4.4 | 15 | 22.3 | 1.2 | 21.2 | 0.00 | 1.1 | 15 | 813 |
| 60 | 1971.1 | 28 | 4.0 | 11 | 21.9 | 0.8 | 21.0 | 0.00 | 1.2 | 11 | 509 |

Note: 1. The baseline scenario. No customers were rejected in all cases. 2. The CPU time for a 10-minute frequency is 34.8 hours.

### 6.7.2 Number of transit lines

A denser transit network provides higher accessibility for customers, but also makes the EIDARP more complex. We evaluate service performance with varying numbers of transit lines (representing network density) within a 16 km by 16 km operational area. Figure 8 illustrates the layouts of four different transit networks with the number of lines ranging from 1 to 4. For transit networks with more than two lines, the number of transit stations is 3 for each transit line, with one transfer station positioned in the middle. Each transit line operates in both directions with a frequency of 30 minutes during a two-hour operational period.

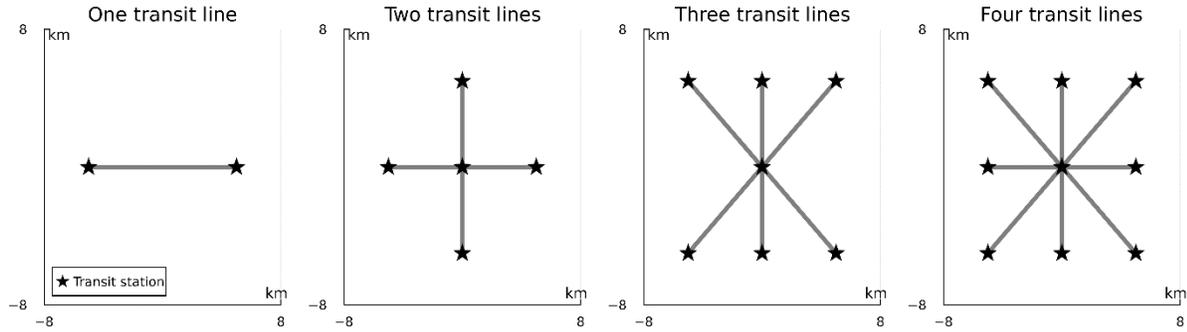

**Figure 8:** The layouts of four transit networks with an increasing number of transit lines.

Table 16 presents the results with the four transit networks and the no-transit scenario. The latter represents a DARP service where buses serve all customers exclusively, resulting in the highest BTT. As the number of transit lines increases, both BTT and CTT decrease. This is because more customers use transit services, with the number increasing from 0 to 58. This also leads to an increase in CTT-transit and CTT-walk. While more customers utilize the transit service with more transit lines, the number of buses used does not necessarily decrease, as customers are served by different buses for first-mile and/or last-mile services. Moreover, the average number of served customers per bus increases from 3.6 to 5.2, suggesting more customers can share bus rides at their first and last miles. The complexity of the EIDARP also depends on the transit network density, reflected in the higher CPU time as the number of transit lines increases. For the no-transit scenario, the CPU time to solve the EIDARP is significantly smaller compared with other cases.

**Table 16:** The impacts of the density of the transit network on the performance of the EIDARP.

| # transit line | BTT | # used buses | # cus/bus | # cus-transit | CTT | CTT-transit | CTT-bus | CTT-walk | WT | # TO4 | CPU (s) |
|---|---|---|---|---|---|---|---|---|---|---|---|
| None | 1971.1 | 28 | 3.6 | 0 | 21.8 | 0.0 | 21.8 | 0.00 | 0.9 | 0 | 148 |
| One | 1930.9 | 27 | 3.9 | 5 | 22.0 | 0.5 | 21.4 | 0.07 | 1.0 | 4 | 625 |
| Two[1] | 1902.9 | 29 | 4.6 | 39 | 21.4 | 2.8 | 18.4 | 0.15 | 2.7 | 36 | 2781 |
| Three | 1819.8 | 31 | 4.7 | 56 | 19.2 | 4.0 | 14.9 | 0.33 | 3.2 | 51 | 10882 |
| Four | 1650.4 | 28 | 5.2 | 58 | 18.9 | 5.1 | 13.1 | 0.72 | 1.8 | 50 | 40663 |

Note: 1. The baseline scenario. No customers were rejected, and no recharging operations were required for all the cases.

# 7     Conclusion

The IDRT service has received increasing attention in recent years. This innovative transportation service helps DRT service operators reduce unnecessary VKT and increase ridership for existing mass transit services. This paper considers an EIDARP, focusing on transit synchronization and the utilization of EVs. A hybrid LNS metaheuristic is developed to efficiently solve the problem for larger problem sizes. The proposed algorithm operates on a departure-expanded network and identifies the feasible origin-destination transit pairs for customers, thereby enhancing the algorithm's efficiency. A nine-step solution feasibility evaluation scheme ensures customers' maximum travel time, transit departure synchronization, and maximum waiting time at transit stations. This scheme calculates the forward slack time to identify candidate insertion positions for vehicle charging operations. The hybrid LNS incorporates a set of local search operators and a DA-based acceptance criterion to balance the diversification and intensification in the metaheuristic design. Several problem-specific LNS operators are developed to efficiently solve the EIDARP. We compared the performance of the hybrid LNS with the exact solution obtained by the state-of-the-art solver (Gurobi) using a set of test instances. Due to the complexity of the EIDARP, the exact (best feasible) solutions can be obtained for the instances with up to 50 customers, given an 8-hour computation time limit. The results show that the hybrid LNS significantly outperforms the solutions found by Gurobi with a short computation time (a couple of minutes). We further tested the hybrid LNS on larger test instances with up to 100 customers, given low and full initial SoCs of vehicles. The results show that the proposed metaheuristic can solve these instances within a reasonable computation time limit.

    We conducted a series of experiments to evaluate the EIDARP system with respect to several key system parameters. The results provide insights for operators to navigate the trade-off between operational costs and customer inconvenience. The findings suggest that the weights in the objective function (to trade off the operational costs and customer inconvenience) and the detour factor significantly influence the outcomes. The results show that using a higher weight on customer travel time can achieve similar system performance compared to using a lower detour factor. Although a lower bus VKT can be achieved by increasing the detour factors (or reducing the weight on customer travel time), it negatively impacts the customer experience. While the integrated service reduces bus VKT compared to the non-integrated service (e.g., DARP), the number of buses (vehicles) used is not necessarily reduced because different buses must serve customers' first and last miles. However, it was found that using a relatively small fleet size would lead to a greater customer rejection rate, indicating that a suitable fleet size needs to be determined to minimize overall system costs. Imposing a maximum waiting time at transit stations does not appear to be critical for improving the system performance, as the customer's maximum travel time constraint (detour) would minimize customer waiting time at transit stations.

    For EIDARP operators, selecting an appropriate service area is crucial for the efficient operation of an EIDARP service. One key consideration is that the bus operational speed should be lower than the transit operational speed. If the bus operational speed is similar to the transit operational speed, fewer customers would likely switch from DRT to integrated DRT due to additional transfer times at transit stations. Additionally, the EIDARP service operates more efficiently in areas with high transit connectivity (accessibility). Although higher transit frequencies benefit both customers and operators, these gains are smaller than those obtained by expanding transit coverage (i.e., adding more lines). Moreover, increased frequency incurs substantial computational costs. Our results suggest that a fully charged bus fleet at the beginning of service ensures higher vehicle availability and reduces the number of buses in use. When scheduling charging operations on scarce charging infrastructure, higher computational time is required to meet charging station capacity constraints.

However, the proposed algorithm still faces challenges when scaling to large real-world problems. When the charging demand is high and the charging infrastructure is limited, charging scheduling becomes a major bottleneck due to the rigid charging capacity constraint. In practice, one can relax this assumption by allowing vehicles waiting at charging stations to adjust their charging schedules to reduce the difficulty of finding feasible schedules. Developing more efficient heuristics using parallel-execution schemes could be another avenue of research. In our numerical studies, more than 50% of the computational time is spent on the repair operators for the instance with 100 requests, suggesting that significant speedups can be achieved for these operators by introducing parameters to control the trade-off between solution quality and computational time. Systematic evaluation of the EIDARP services for a real-world case study could be another interesting research area to explore. Furthermore, the impact of the maximum transfer time threshold on the computational efficiency of the hybrid metaheuristic for larger transit networks needs to be investigated to develop more efficient algorithms.

Several future research directions can be explored by extending the current study. First, customer inconvenience can be further improved by considering the constraints of customer transfer times and the number of transfers within the transit network. Second, a more realistic non-linear charging function can be incorporated for charging schedule modeling. Jointly optimizing charging station location and demand-responsive bus routing can further reduce system costs and increase service quality. Extending the static EIDARP to a dynamic one with the development of efficient solution algorithms would be a relevant research area for its practical applications. Finally, integrating different sources of uncertainty (e.g., travel time, energy consumption, and transit delays) in the modeling approach could further improve the reliability and robustness of the integrated DRT service.

**Acknowledgments:**

The work was supported by the Luxembourg National Research Fund (C20/SC/14703944).

## Appendix A.    Notation

| Sets | |
|---|---|
| $K$ | Set of electric buses |
| $O, \bar{O}$ | Set of origin and destination depots |
| $P$ | Set of customers' origins, i.e. $P = \{1, \ldots, n\}$ |
| $R$ | Set of customers, i.e., $R = \{1, \ldots, n\}$ |
| $D$ | Set of customers' destinations, i.e., $D = \{n+1, \ldots, 2n\}$ |
| $\mathcal{L}$ | Set of transit lines |
| $\mathcal{F}$ | Set of physical transit stations, i.e., $F = \{F_l\}_{l \in \mathcal{L}}$ |
| $\mathcal{D}$ | Set of transit (train) departures, i.e., $\mathcal{D} = \{\mathcal{D}_l\}_{l \in \mathcal{L}}$ |
| $G$ | Set of transit (dummy) nodes, i.e., $G = \{g_d^f\}_{f \in F_l, d \in \mathcal{D}_l, l \in \mathcal{L}}$ |
| $C$ | Set of charging stations |
| $S$ | Set of dummy charger nodes, i.e., $S = \{S_c\}_{c \in C}$ |
| $\mathcal{T}_r$ | Set of transit pairs of customer $r$ |
| $\mathcal{A}_B$ | Set of bus arcs |
| $\mathcal{A}_G$ | Set of transit arcs, each of which represents a transit pair from one station to another, as generated by Algorithm 4. |
| $\mathcal{A}_R$ | Set of walking arcs between origins/destinations and transit nodes within a constant maximum walking distance |
| $N$ | Set of all nodes. $N = P \cup D \cup G \cup S \cup \{O, \bar{O}\}$ |
| $\mathcal{A}$ | Set of all arcs. $\mathcal{A} = \mathcal{A}_B \cup \mathcal{A}_G \cup \mathcal{A}_R$ |
| **Parameters** | |
| $n$ | Number of customers |
| $t_{ij}$ | Travel time from bus node $i$ to bus node $j$, $\forall (i,j) \in \mathcal{A}_B$ |
| $t_{ij}^w$ | Walking time from node $i$ to node $j$, $\forall (i,j) \in \mathcal{A}_R$ |
| $t_{ij}^m$ | Travel time from transit node $i$ to $j$, $\forall (i,j) \in \mathcal{A}_G$ |
| $t_{end}$ | The end of the service operation duration |
| $c_{ij}$ | Distance from node $i$ to node $j$ |
| $\bar{c}$ | The average distance of the arcs in $\mathcal{A}$ |
| $c_{max}$ | The distance of the longest arc in $\mathcal{A}$ |
| $u_i$ | Service time at node $i$ |
| $e_i, l_i$ | Earliest and latest starting times of service at bus node $i$ |
| $\bar{\theta}_i, \underline{\theta}_i$ | Arrival and departure time of transit vehicle (train) at node $i$, $i \in G$ |
| $\gamma$ | Maximum waiting time for inter-modal transfer |
| $\eta^{min}, \eta^{max}$ | The minimum and maximum transfer time thresholds to trim off infeasible transfer arcs within the same transit station in a transit network |
| $\varphi$ | Detour factor |
| $Q^k$ | Capacity of bus $k$ |
| $\beta^k$ | Energy consumption rate per kilometer traveled for bus $k$ |
| $\alpha_s$ | Charging speed of charger $s$ |
| $E_{min}^k, E_{max}^k, E_{init}^k$ | Minimum, maximum, and initial SoC of bus $k$ |
| $L_r^{max}$ | Maximum travel time for customer $r$. $L_r^{max} = t_{r,n+r} \times \varphi$ |
| $t_{end}$ | End of service operation duration |
| $\lambda_1, \lambda_2$ | Weighted coefficients in the objective function |
| $\omega$ | Penalty associated with one customer rejection |
| $T_{max}$ | The maximum temperature for the deterministic annealing. $T_{max} = t_{max} \times \bar{c}$, where $t_{max}$ is a user-defined parameter. |
| $T_{red}$ | Cooling temperature for the deterministic annealing |
| $M$ | Large positive number |
| **Auxiliary variables** | |
| $A_i^k$ | Arrival time of bus $k$ at a transit node $i$ |
| $E_i^k$ | Battery energy level of bus $k$ when arriving at node $i$ |
| $h_{kk'}^{ss'}$ | 1 if bus $k'$ arrives at dummy charging node $s'$ later than bus $k$'s arrival time at dummy charging node $s$, 0 otherwise |
| $dep_r, arr_r$ | Departure and arrival time of customer $r$, respectively |
| $v_r$ | 1 if customer $r$ is served, 0 otherwise |
| **Decision variables** | |
| $x_{ij}^k$ | 1 if arc $(i,j)$ is traversed by bus $k$, 0 otherwise. |
| $y_{ij}^r$ | If a bus serves customer $r$ on arc $(i,j)$, 0 otherwise |
| $z_{ij}^r$ | 1 if customer $r$ uses transit arc $(i,j)$, 0 otherwise |
| $w_{ij}$ | 1 if customer walks from node $i$ to node $j$, 0 otherwise |
| $B_i^k$ | Beginning time of service of bus $k$ at node $i$ |
| $\tau_s^k$ | Charging duration for bus $k$ at charger $s \in S$ |

## Appendix B. MILP formulation of the EIDARP

The EIDARP is formulated as a MILP as follows. We use a three-index formulation to track buses' routes and energy consumption as we consider a fleet of heterogeneous buses. Let $\sigma_B^-(i)$ and $\sigma_B^+(i)$ denote the set of incoming and outgoing nodes of node $i$ in bus arc set $\mathcal{A}_B$. Similarly, $\sigma_G^-(i)$ and $\sigma_G^+(i)$ refer to the incoming and outgoing nodes of node $i$ in the transit arc set $\mathcal{A}_G$, and $\sigma_R^-(i)$ and $\sigma_R^+(i)$ those in the walking arc set $\mathcal{A}_R$. The decision variables include: whether a bus travels from $i$ to $j$ ($x_{ij}^k$), a customer travels from $i$ to $j$ ($y_{ij}^r$); walks from $i$ to $j$ ($w_{ij}$); whether a leg of transit line is used for a customer from $i$ to $j$ ($z_{ij}^r$); bus charging time at a charger $s$ ($\tau_s^k$). Given a set of requests, the objective function (C.1) minimizes the weighted sum of bus routing time (first term), customers' in-vehicle (including bus and transit vehicles) travel time and walking time (second term), and the penalty of unserved customers (third term).

$$Min \; \lambda_1 \sum_{(i,j) \in A_B} \sum_{k \in K} t_{ij} x_{ij}^k + \lambda_2 \sum_{r \in R} \bar{L}_r + \omega \sum_{r \in R} (1 - v_r) \qquad (C.1)$$

where $\bar{L}_r$ is the summation of customer $r$'s in-bus time, in-transit (train) time, and walking time, calculated as:

$$\bar{L}_r = \sum_{(i,j) \in A_B} t_{ij} y_{ij}^r + \sum_{(i,j) \in A_G} t_{ij}^T z_{ij}^r + \sum_{i \in G} t_{ri}^w w_{ri} + \sum_{i \in G} t_{n+r,i}^w w_{n+r,i} \; \forall r \in R \qquad (C.2)$$

subject to

$$\sum_{j \in \sigma_B^+(o^k)} x_{o^k j}^k = 1, \forall k \in K \tag{C.3}$$

$$\sum_{i \in \sigma_B^-(\bar{o}^k)} x_{i,\bar{o}^k}^k = 1, \forall k \in K \tag{C.4}$$

$$\sum_{j \in \sigma_B^+(o)} x_{oj}^k = 0, \forall k \in K, o \in O, o \neq o^k \tag{C.5}$$

$$\sum_{j \in \sigma_B^-(i)} x_{ji}^k - \sum_{j \in \sigma_B^+(i)} x_{ij}^k = 0, \forall k \in K, i \in N \setminus \{O, \bar{O}\} \tag{C.6}$$

$$y_{ij}^r \leq \sum_{k \in K} x_{ij}^k, \forall r \in R, (i,j) \in \mathcal{A}_B \tag{C.7}$$

Constraints (C.3) - (C.4) ensure that buses need to leave and return to their respective depots ($o^k$). Constraints (C.5) forbid each bus to visit other depots except its own. Constraints (C.6) ensure bus flow conservation. Constraints (C.7) ensure the consistency between $x_{ij}^k$ and $y_{ij}^r$.

$$w_{ij} = 0, \forall i \in P, j \in G, (i,j) \notin \mathcal{A}_R \tag{C.8}$$

$$v_i = \sum_{j \in \sigma_R^+(i)} w_{ij} + \sum_{j \in \sigma_B^+(i)} y_{ij}^i, \forall i \in P \tag{C.9}$$

$$v_i = \sum_{j \in \sigma_R^-(n+i)} w_{j,n+i} + \sum_{j \in \sigma_B^-(n+i)} y_{j,n+i}^i, \forall i \in P \tag{C.10}$$

$$\sum_{j \in \sigma_B^+(i)} y_{ij}^r - \sum_{j \in \sigma_B^-(i)} y_{ji}^r = 0, \forall i \in P \cup D, r \in R, i \neq r, i \neq n+r \tag{C.11}$$

$$w_{i,n+r} - w_{ri} + \sum_{j \in \sigma_B^+(i)} y_{ij}^r - \sum_{j \in \sigma_B^-(i)} y_{ji}^r + \sum_{j \in \sigma_G^+(i)} z_{ij}^r - \sum_{j \in \sigma_G^-(i)} z_{ji}^r = 0, \forall r \in R, i \in G \tag{C.12}$$

Constraints (C.8) enforce that customers cannot walk from/to transit stations if the distance between their origins/destinations and transit stations exceeds the maximum walking distance threshold. $\mathcal{A}_R$ is the set of directed walking arcs to connect transit station nodes within the maximum walking distance threshold. Constraints (C.9) - (C.10) state that each customer can be either served or rejected. If served ($v_i = 1$), the customer must be served by a bus or walk to/from a transit station if their origin and/or destination is located within the maximum walking distance threshold. Note that given customers' time windows constraints at origins/destinations, their maximum ride times, and the inter-modal allowable transfer buffer (i.e., 15 minutes before/after the departure/arrival of trains) associated with transit nodes, a preprocessing procedure is applied to trim off infeasible arcs by verifying these time-window constraints. Consequently, $\sigma_B^+(i)$ and $\sigma_B^-(i)$ contain only relevant (feasible) outgoing/incoming bus arcs. Customer flow consistency is ensured at bus nodes (Eq. C.11) and transit nodes (Eq. C.12).

$$\sum_{(i,j) \in A_G} z_{ij}^r \leq 1, \forall r \in R \tag{C.13}$$

$$w_{ri} \leq \sum_{j \in \sigma_G^+(i)} z_{ij}^r, \forall r \in R, i \in G \tag{C.14}$$

$$\sum_{j \in \sigma_B^+(i)} x_{ij}^k \leq \sum_{j \in \sigma_B^+(n+i)} x_{n+i,j}^k + M_1 \sum_{(i,j) \in A_G} z_{ij}^i, \forall k \in K, i \in P \tag{C.15}$$

Constraints (C.13) state that each customer can use a timetabled transit service at most once with entry transit node $i$ and exit transit node $j$. If customers walk to a transit station, they must be served by transit services (Eq. C.14). Constraints (C.15) state that if customers do not use a transit service, they must be served by the same buses that visit their origins and destinations. Note that we use a set of big-Ms (sufficiently large positive numbers) to model the if-then conditions.

$$y_{ij}^r = 0, \forall r \in R, i \in S \cup O, j \in \sigma_B^+(i) \tag{C.16}$$

$$\sum_{r \in P} y_{ij}^r \leq Q^k + M_2(1 - x_{ij}^k), \forall k \in K, (i,j) \in \mathcal{A}_B \tag{C.17}$$

$$B_j^k \geq B_i^k + t_{ij} + u_i - M_3(1 - x_{ij}^k), \forall k \in K, (i,j) \in \mathcal{A}_B, i \notin S \tag{C.18}$$

$$B_j^k \geq B_i^k + t_{ij} + \tau_i^k + u_i - M_3(1 - x_{ij}^k), \forall k \in K, (i,j) \in \mathcal{A}_B, i \in S \tag{C.19}$$

$$e_i \leq B_i^k \leq l_i, \forall k \in K, i \in P \cup D \cup \bar{O} \tag{C.20}$$

Constraints (C.16) ensure that no customer can be on board at the depot and at the chargers, while constraints (C.17) state that the passenger load cannot exceed the bus capacity ($Q^k$). Constraints (C.18) - (C.19) ensure consistency of the beginning time of service ($B_i^k$) at bus nodes and charger nodes, respectively. Constraints (C.20) state the time-window constraints at pickup nodes and drop-off nodes.

$$\underline{\theta}_i - \gamma \leq A_i^k \leq \bar{\theta}_i + \gamma, \forall k \in K, i \in G \tag{C.21}$$

$$A_j^k \geq B_i^k + t_{ij} + u_i - M_3(1 - x_{ij}^k), \forall k \in K, j \in G, i \in \sigma_B^-(j) \tag{C.22}$$

$$A_j^k \leq B_i^k + t_{ij} + u_i + M_3(1 - x_{ij}^k), \forall k \in K, j \in G, i \in \sigma_B^-(j) \tag{C.23}$$

$$A_i^k \geq \underline{\theta}_i - \gamma - M_3\left(1 - \sum_{j \in \sigma_G^+(i)} z_{ij}^r\right), \forall r \in R, i \in G, k \in K \tag{C.24}$$

$$A_i^k \leq \underline{\theta}_i + M_3\left(1 - \sum_{j \in \sigma_G^+(i)} z_{ij}^r\right), \forall r \in R, i \in G, k \in K \tag{C.25}$$

$$B_i^k \geq \bar{\theta}_i - M_3\left(1 - \sum_{j \in \sigma_G^-(i)} z_{ji}^r\right), \forall r \in R, i \in G, k \in K \tag{C.26}$$

$$B_i^k \leq \bar{\theta}_i + \gamma + M_3\left(1 - \sum_{j \in \sigma_G^-(i)} z_{ji}^r\right), \forall r \in R, i \in G, k \in K \tag{C.27}$$

As buses can visit a transit station for dropping off and/or picking up customers, constraints (C.21) define the time windows associated with transit nodes as the departure time of a transit service minus and plus a maximum waiting time. To ensure a maximum transfer (waiting) time at transit stations, an auxiliary variable tracking bus arrival time from node $i$ to $j$ is introduced by constraints (C.22) and (C.23). Constraints (C.24) - (C.27) are inter-modal transfer time constraints at transit stations for customers transfer between buses and transit services. Constraints (C.24) - (C.25) coordinate a bus's arrival at transit node $i$ for the first-mile service within the predefined time window $[\underline{\theta}_i - \gamma]$ with a maximum transfer time $\gamma$. While constraints (C.26) - (C.27) are used for the last-mile service to ensure that buses pick up customers from transit nodes (stations) within the buffer time $\gamma$ after a transit vehicle arrival at time $\bar{\theta}_i$.

$$\underline{\theta}_j \leq l_i + t_{ij}^w + M_3(1 - w_{ij}), \ \forall i \in P, j \in \sigma_B^+(i) \tag{C.28}$$

$$\underline{\theta}_j \geq e_i + t_{ij}^w - M_3(1 - w_{ij}), \ \forall i \in P, j \in \sigma_B^+(i) \tag{C.29}$$

$$\bar{\theta}_j \leq l_i - t_{ji}^w + M_3(1 - w_{ji}), \ \forall i \in D, j \in \sigma_B^-(i) \tag{C.30}$$

$$\bar{\theta}_j \geq e_i - t_{ji}^w - M_3(1 - w_{ji}), \ \forall i \in D, j \in \sigma_B^-(i) \tag{C.31}$$

Constraints (C.28) - (C.31) ensure time-window constraints if customers walk from their origins to transit stations or from transit stations to their destinations.

$$B_i^k \leq M_3 \sum_{j \in \sigma_B^+(i)} x_{ij}^k, \ \forall i \in P \cup D, k \in K \tag{C.32}$$

$$dep_r = \sum_{k \in K} B_r^k + \sum_{j \in \sigma_R^+(r)} w_{rj}(\underline{\theta}_j - t_{rj}^w), \ \forall r \in R \tag{C.33}$$

$$arr_r = \sum_{k \in K} B_{n+r}^k + \sum_{j \in \sigma_R^-(n+r)} w_{j,n+r}(\bar{\theta}_j + t_{j,n+r}^w), \ \forall r \in R \tag{C.34}$$

$$arr_r - dep_r \leq L_r^{\max}, \ \forall r \in R \tag{C.35}$$

Constraints (C.32) ensures that the beginning time of service at a customer's origin or destination is set to zero if a bus $k$ does not visit it. Constraints (C.33) and (C.34) calculate customers' departure and arrival times when they are served by buses or walk to/from transit stations. Constraints (35) determines the journey time of customers (i.e., $L_r = arr_r - dep_r$) and ensures that the maximum journey time of customers is not violated. Note that the customers' maximum journey time is determined as the multiplication of a detour factor and the direct travel times by bus between their origins and destinations.

$$E_0^k = E_{\text{init}}^k, \ \forall k \in K \tag{C.36}$$

$$E_{\min}^k \leq E_i^k \leq E_{\max}^k, \ \forall i \in N \setminus O, k \in K \tag{C.37}$$

$$E_j^k \geq E_i^k + \alpha_i \tau_i^k - c_{ij}\beta^k - M_4(1 - x_{ij}^k), \ \forall k \in K, (i,j) \in \mathcal{A}_B, i \in S \tag{C.38}$$

$$E_j^k \leq E_i^k + \alpha_i \tau_i^k - c_{ij}\beta^k + M_4(1 - x_{ij}^k), \ \forall k \in K, (i,j) \in \mathcal{A}_B, i \in S \tag{C.39}$$

$$E_j^k \geq E_i^k - c_{ij}\beta^k - M_4(1 - x_{ij}^k), \ \forall k \in K, (i,j) \in \mathcal{A}_B, i \notin S \tag{C.40}$$

$$E_j^k \leq E_i^k - c_{ij}\beta^k + M_4(1 - x_{ij}^k), \ \forall k \in K, (i,j) \in \mathcal{A}_B, i \notin S \tag{C.41}$$

Constraints (C.36) set the initial state of charge (SoC) of buses. Constraints (C.37) ensure the upper bound and lower bound of buses' SoC. Constraints (C.38) and (C.39) track buses' SoC changes at charger nodes, while constraints (C.40) and (C.41) track buses' SoC changes at bus nodes.

$$\sum_{j \in \sigma_B^+(s)} x_{sj}^k \leq 1, \ \forall k \in K, s \in S \tag{C.42}$$

$$\sum_{j \in \sigma_B^+(l)} x_{lj}^k \leq \sum_{j \in \sigma_B^+(h)} x_{hj}^k, \ \forall k \in K, h, l \in S_c, c \in C, h < l \tag{C.43}$$

$$h_{kk'}^{ss'} + h_{k'k}^{s's} = 1, \forall s, s' \in S_c, c \in C, k, k' \in K, k \neq k' \tag{C.44}$$

$$B_{s'}^{k'} \geq B_s^k + \tau_s^k - M_3(1 - h_{kk'}^{ss'}), \forall s, s' \in S_c, c \in C, k, k' \in K, k \neq k' \tag{C.45}$$

Constraints (C.42) – (C.45) are related to bus charging scheduling constraints to ensure no overlaps of charging activities on each charger. Constraints (C.42) state that each vehicle $k$ can connect each dummy charger node $s$ at most once, but a dummy charging node $s$ can be visited by different buses. With multiple charger dummy nodes associated with a charger, a bus is allowed to charge multiple times at the same charger by visiting different charger dummy nodes. Note that the charging dummy nodes are an ordered set. The symmetry issue is avoided by constraints (C.43), which ensure that the dummy charger nodes used need to follow their ordered indexes. $h_{kk'}^{ss'}$ is an indicator being 1 if bus $k'$ arrives later than bus $k$ at the same physical charger. $s'$ and $s$ are the visited dummy charging nodes for bus $k'$ and $k$, respectively (Eq. C.44). Constraints (C.45) ensure that if bus $k'$ arrives later than bus $k$, bus $k'$ can start its charging session after the end of the charging session of bus $k$.

$$v_r = 1, \ \forall r \in R \tag{C.46}$$

Note that if customer rejection is not allowed, constraints (C.46) can be added.
Constraints (C.47) – (C.55) define the domain of all decision variables and auxiliary variables.

$$x_{ij}^k \in \{0,1\}, \forall k \in K, (i,j) \in \mathcal{A}_B \tag{C.47}$$
$$y_{ij}^r \in \{0,1\}, \forall r \in R, (i,j) \in \mathcal{A}_B \tag{C.48}$$
$$z_{ij}^r \in \{0,1\}, \forall r \in R, (i,j) \in \mathcal{A}_G \tag{C.49}$$
$$w_{ij} \in \{0,1\}, \forall (i,j) \in \mathcal{A}_R \tag{C.50}$$
$$B_i^k \geq 0, E_i^k \geq 0, \forall k \in K, i \in P \cup D \cup G \cup S \cup \{O, \bar{O}\} \tag{C.51}$$
$$\tau_s^k \geq 0, \forall k \in K, s \in S \tag{C.52}$$
$$A_i^k \geq 0, \forall k \in K, i \in G \tag{C.53}$$
$$h_{kk'}^{ss'} \in \{0,1\}, \forall s, s' \in S_c, c \in C, k, k' \in K, k \neq k' \tag{C.54}$$
$$arr_r \geq 0, dep_r \geq 0, v_r \in \{0,1\}, \forall r \in R \tag{C.55}$$

We set tight values for big-Ms as $M_1 = max\{|\sigma_B^+(i)|: i \in P\}, M_2 = n, M_3 = t_{end}, M_4 = max\{E_{max}^k: k \in K\}$.

## Appendix C.  An illustrative example of the departure-expanded graph

We present an illustrative example to explain the departure-expanded graph in Figure C1 and Figure C2. The departure-expanded graph is a variant of the time-expanded graph to model the transit-integrated routing problem. Figure C1 presents a scheduled transit service with 6 physical transit stations and 2 transit lines. Line 1 (A-B-C-D) and Line 2 (D-F-G-H) are connected at two transfer stations (B/F and D/H). Transit line nodes B and F (D and H) have the same physical transit stations (stops). The scheduled timetables for each line are presented on the left of Figure C1. Compact timetables containing both directions are shown on the left (usually separated into two timetables, one for each direction). Line 1 has three ordered departures (schedules), while Line 2 has two ordered departures. The arrows indicate the direction, and the numbers in the timetables denote the departure times of the transit vehicle at each station. Figure C2 illustrates the structure of the corresponding departure-expanded graph with the corresponding transit nodes $G = \{1,2,\ldots,20\}$. An arc means that there is a direct route without transfer from one station to another of the same departure. Given the timetables of transit lines, transit stations visited by a transit line are duplicated according to the number of departures. Ordered transit nodes are structured as a departure-expanded graph in which direct arcs connect two transit nodes when they can be served by the same vehicle (same transit line) or when the allowable transfer time between different transit lines is within a predefined transfer time limit related to the quality of service of the operator. The numbers within filled points are ordered indices of transit nodes. Each transit node is associated with its respective (transit vehicle) departure time defined by the timetable in Figure C1.

For example, the direct arc from node 1 to node 3 in Figure C2 means that a transit vehicle traverses consecutively stations 1, 2, and 3. The weights on the arcs correspond to the travel times between two connected stations. There are 6 direct arcs for the first departure of each line. The generation of transfer arcs at transfer stations are explained in Figure C2 (b). For the first transfer station, three transit nodes (2, 10, and 18) are associated with the station node B (line 1), and two transit nodes (6 and 14) are associated with the station node F (line 2). The transfer from node 2 (departure time is 7:26) to node 6 (departure time is 7:29) takes 3 minutes, while the transfer from node 2 (departure time is 7:26) to node 14 (departure time is 7:40) takes 14 minutes. By limiting the maximum allowed transfer time (e.g., $\eta^{max}$=10 minutes) within the same transit station, infeasible transfer arcs are trimmed off. As a result, there are only two feasible transfer arcs are created given the 10-minute maximum transfer time limit (Figure C2).

Algorithm 4 describes the details to generate the departure-expanded graph $\mathcal{G}_T$. The algorithm includes mainly two steps: the first step is to generate a temporary departure-expanded transit graph $\mathcal{G}'_T = (G, \mathcal{A}'_G)$ with temporary transit pairs (TPs) $\mathcal{A}'_G$ connecting direct transit and feasible transfer arcs for feasible shortest path travel time calculation between pairs of transit nodes. The inputs are the sets of transit nodes $G$, transit lines $\mathcal{L}$, departures $\mathcal{D}$, transit arrival times $\bar{\theta}$, and departure times $\underline{\theta}$ at transit nodes given the timetables, and minimum and maximum allowable transfer times $\eta^{min}$ and $\eta^{max}$, respectively. The output is a departure-expanded graph and the feasible shortest travel times on TPs, i.e., $t^m_{ij}, \forall i,j \in G$. For the transit line and departure pair, the direct arcs and transfer arcs are inserted into $\mathcal{A}'_G$. Then we apply the conventional Dijkstra Algorithm (Dijkstra, 1959) to get the shortest travel time between pairs of transit nodes on $\mathcal{G}'_T$. Finally, transfer arcs are removed.

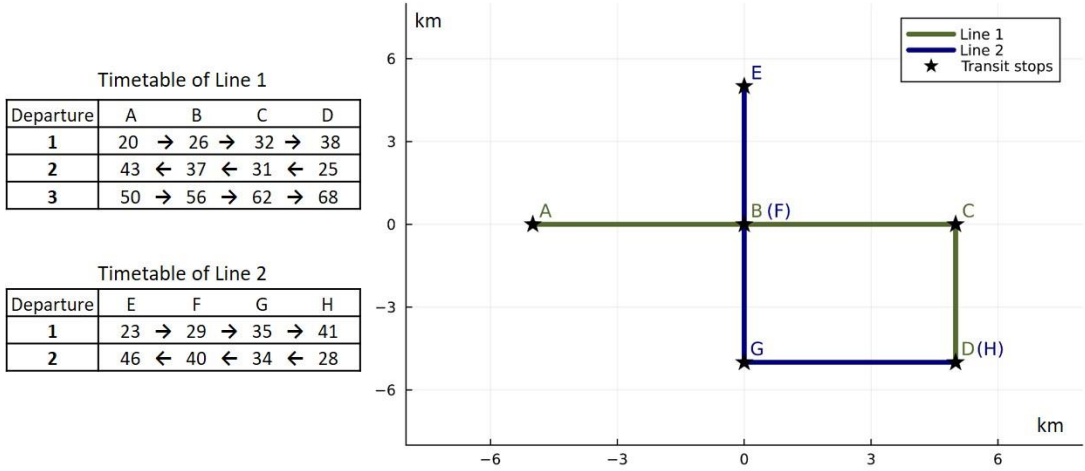

**Figure C1:** An illustrative example of transit service timetables with two transit lines: Line 1 (A-B-C-D) and Line 2 (E-F-G-H).

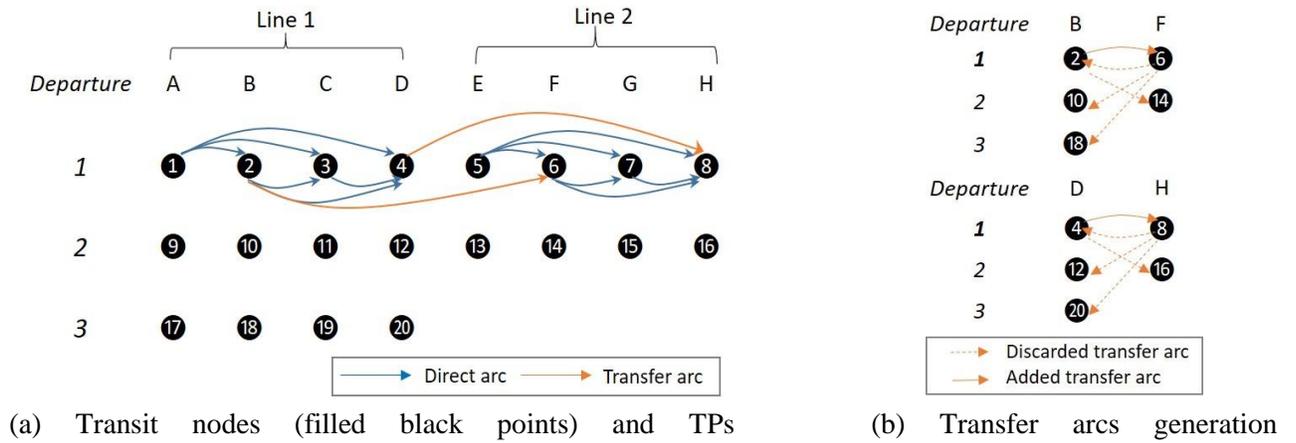

(a) Transit nodes (filled black points) and TPs     (b) Transfer arcs generation

**Figure C2:** Transit arc generation for the departure-expanded graph of the example in C1.

The density of the transit network has a pronounced effect on the complexity of the EIDARP. For instance, consider adding a second transit line $l_2$ (with $|\mathcal{D}_{l_2}|$ departures for $|\mathcal{F}_{l_2}|$ transit stations) to an existing single-line transit network with $l_1$ (with $|\mathcal{D}_{l_1}|$ departures for $|\mathcal{F}_{l_1}|$ transit stations). This expansion introduces:

1) $|\mathcal{D}_{l_2}||\mathcal{F}_{l_2}|$ transit nodes,
2) Intra-line TPs for $l_2$: $|\mathcal{D}_{l_2}|(|\mathcal{F}_{l_2}|-1)$,
3) Inter-line TPs (between $l_1$ and $l_2$): $2|\mathcal{D}_{l_2}||\mathcal{F}_{l_2}| \times |\mathcal{D}_{l_1}||\mathcal{F}_{l_1}|$ in the worst case (two directions for each pair of transit node combinations).

Algorithm 4 reduces the inter-line TPs by pruning infeasible transfers based on the minimum and maximum transfer times and taking only the shortest path between every two transit nodes. Moreover, introducing $l_2$ also induces additional bus arcs connecting depots, customer origins and destinations, and charging stations. For $n$ customers and $|S|$ charging nodes, it generates $2(2n+2+|S|) \times |\mathcal{F}_{l_2}| \times |\mathcal{D}_{l_2}|$ additional bus arcs.

| | Algorithm 4. Transit pairs generation. |
|---|---|
| | **Input:** $G, \mathcal{L}, \mathcal{D}, \bar{\theta}, \underline{\theta}, \eta^{min}, \eta^{max}$ |
| | **Output:** $\mathcal{A}_G, t_{ij}^m \quad \forall i, j \in G$ |
| 1: | Set $\mathcal{A}_G' = \emptyset, \mathcal{A}_G = \emptyset, t_{ij}^m = \infty, \forall j \in G$ |
| 2: | **for** $l \in \mathcal{L}$ |
| 3: |    **for** $d \in \mathcal{D}_l$ |
| 4: |       Insert direct arcs $(i,j)$ to $\mathcal{A}_G'$ and set $t_{ij}^T = \underline{\theta}_j - \underline{\theta}_i$ |
| 5: |       **for** $i \in$ transfer nodes of line $l$ |
| 6: |          **for** $(i,j) \in$ outgoing transfer arcs from $i$ |
| 7: |             **if** $\eta^{min} \leq \underline{\theta}_j - \bar{\theta}_i \leq \eta^{max}$ |
| 8: |                Insert transfer arcs $(i,j)$ to $\mathcal{A}_G'$ and set $t_{ij}^T = \underline{\theta}_j - \bar{\theta}_i$ |
| 9: |             **end if** |
| 10: |          **end for** |
| 11: |       **end for** |
| 12: |    **end for** |
| 13: | **end for** |
| 14: | $\mathcal{G}_T' = (G, \mathcal{A}_G')$ |
| 15: | **for** $i \in G$ |
| 16: |    Calculate the shortest travel time $\tilde{t}_{ij}^T$ for $\forall j \in G, j \neq i$ using Dijkstra Algorithm |
| 17: |    Insert $(i,j)$ to $\mathcal{A}_G$ and update $t_{ij}^T = \tilde{t}_{ij}^T$ if $t_{ij}^T < \infty$ |
| 18: | **end for** |
| 19: | Remove transfer arcs from $\mathcal{A}_G$ |

## Appendix D.     An illustrative example of customer journeys and bus paths

Figure D1 illustrates customer journeys and bus paths for the five travel options (TOs) introduced in Section 3.1. Three buses serve these customers. Bus 2 serves the first-mile trip legs for TO2 (bus-transit-walk) and TO4 (bus-transit-bus), while Bus 3 covers the last-mile legs for TO3 (walk-transit-bus) and TO4 (bus-transit-bus). The customer using TO5 (bus only) is served door-to-door by Bus 1.

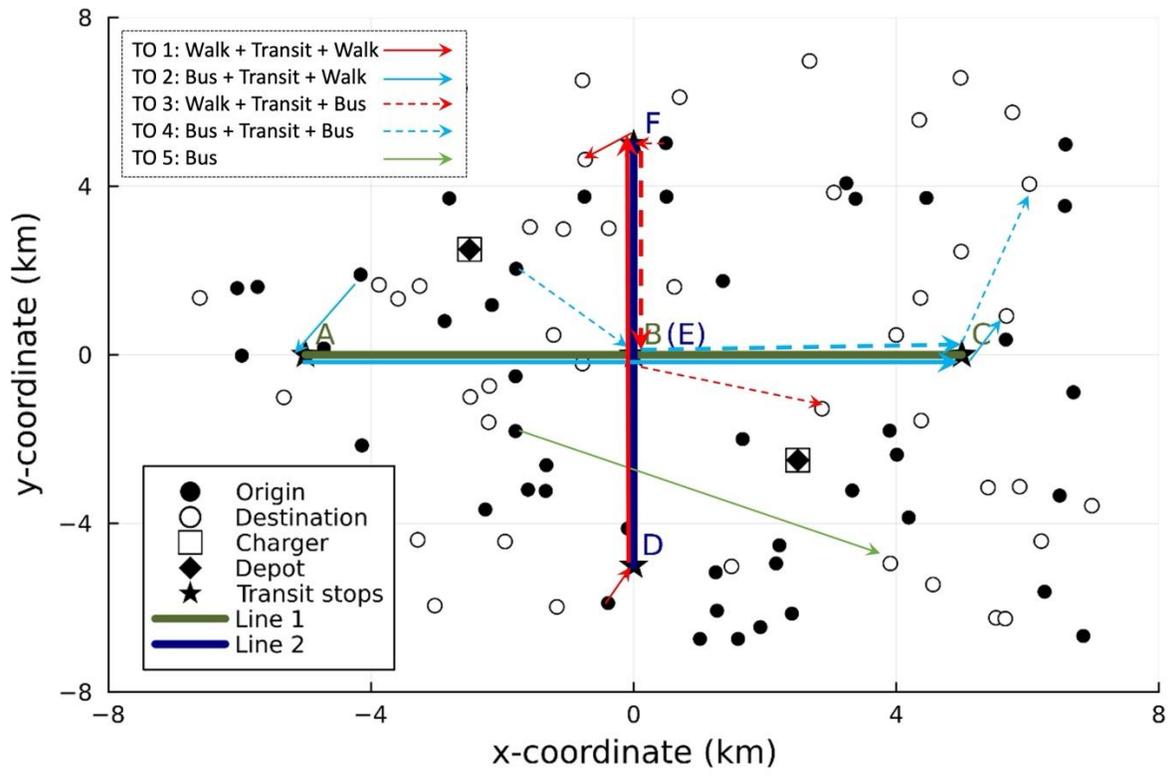

(a) Extract of customer journeys using five different travel options (TOs) in the solution.

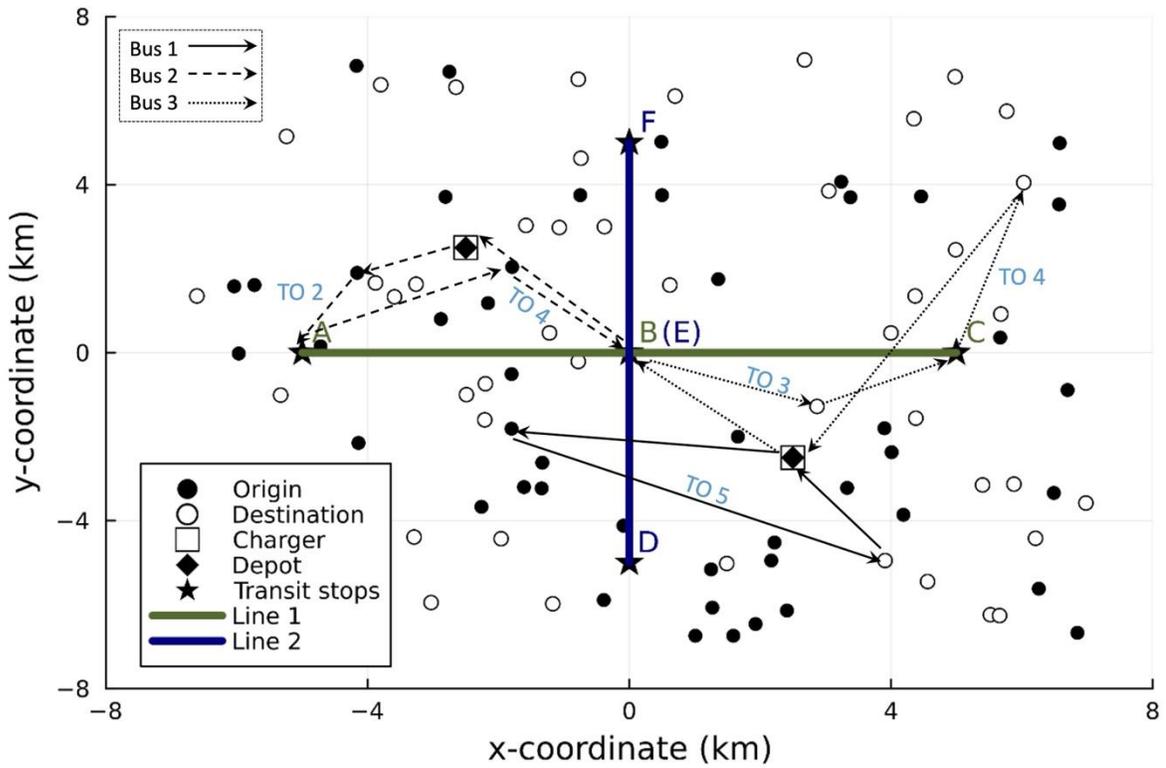

(b) Extract of three different buses' paths in the solution.

**Figure D1:** An illustrative example of customer journeys (a) and bus paths (b).